%% file: HSL_tori_final.tex
\documentclass[11pt]{amsart}
\usepackage[pdfpagemode=None,colorlinks=true,linkcolor=blue,citecolor=blue]{hyperref}
\usepackage{verbatim, amscd,epsfig,amsmath,amsthm,amstext,amssymb,
hyperref 
}
 \usepackage{bbm}%,theorem}
\usepackage{microtype}

\theoremstyle{plain}
\newtheorem*{theorem*}{Theorem}
\newtheorem{theorem}{Theorem}[section]
\newtheorem{lemma}[theorem]{Lemma}
\newtheorem{cor}[theorem]{Corollary}
\newtheorem{prop}[theorem]{Proposition}

\theoremstyle{definition}
\newtheorem{definition}[theorem]{Definition}

\theoremstyle{remark}
\newtheorem{rem}[theorem]{Remark}
\newtheorem{remark}[theorem]{Remark}

\newcommand{\assign}{:=}
\newcommand{\ttrivial}{\widetilde{\trivial{}}}

\newcommand{\tmop}[1]{\ensuremath{\operatorname{#1}}}

\numberwithin{equation}{section}  
\renewcommand{\Re}{{\rm Re}\,}
\renewcommand{\Im}{{\rm Im}\,}

\newcommand{\ejbh}{e^{\frac{j\beta}2}}
\newcommand{\emjbh}{e^{-\frac{j\beta}2}}
\newcommand{\ebnh}{e_{\frac{\beta_0}{2}}}
\newcommand{\embnh}{e_{\frac{-\beta_0}{2}}}
\newcommand{\edeltaB}{e^{2 \pi i\left\langle \delta-B, z \right\rangle}}
\newcommand{\edelta}{e^{2 \pi i\left\langle \delta, z \right\rangle}}
\newcommand{\Cj}{\mathfrak{C}}

\newcommand{\R}{\mathbb{R}}
\newcommand{\C}{\mathbb{ C}}
\newcommand{\Z}{\mathbb{ Z}}

\renewcommand{\H}{\mathbb{ H}}

\newcommand{\N}{\mathbb{ N}}
\renewcommand{\P}{\mathbb{ P}}
\newcommand{\HP}{\H\P}
\newcommand{\CP}{\C\P}
\renewcommand{\L}{{\mathcal L}}
\newcommand{\E}{{\mathcal E}}

\newcommand{\trivial}[1]{\underline{\H}^{#1}}

\newcommand{\invers}{^{-1}}
\DeclareMathOperator{\End}{End}
\DeclareMathOperator{\Eig}{Eig}
\DeclareMathOperator{\Hom}{Hom}

\DeclareMathOperator{\Span}{Span}
\DeclareMathOperator{\Id}{Id}

\DeclareMathOperator{\Spec}{Spec}
\DeclareMathOperator{\Gr}{Gr}

\newcommand{\Hh}{{\mathcal H}}

\begin{document}
\title[Hamiltonian Stationary Lagrangian tori]{Darboux transforms and spectral curves of Hamiltonian Stationary Lagrangian tori}
\date{\today}
\author{K. Leschke, P. Romon}
\address{Katrin Leschke\\ Department of Mathematics\\University of Leicester\\University Road\\
Leicester, LE1 7RH, UK}

\address{Pascal Romon\\Universit\' e Paris-Est Marne-la-Vall\' ee\\
5 boulevard Descartes\\ 
77454 Marne la Vall\' ee cedex 2, FRANCE}

\email{k.leschke@mcs.le.ac.uk, pascal.romon@univ-mlv.fr}
\thanks{First author supported by DFG SPP 1154 ``Global Differential Geometry''}
\maketitle

\input{abstract}

%\input{spectral}

\input{intro}

\input{holomorphic}
\input{darboux}
\input{example}
\input{mudarboux}
\input{spectralcurve}
 
\bibliographystyle{amsplain}
\bibliography{doc}
%erzeugen der bibliography durch bibtex ``documentname''

\end{document}

%% file: abstract.tex
\begin{abstract}
  The multiplier spectral curve of a conformal torus $f: T^2\to S^4$
  in the 4--sphere is essentially \cite{conformal_tori} given by all
  Darboux transforms of $f$.  In the particular case when the
  conformal immersion is a Hamiltonian stationary torus $f: T^2
  \to\R^4$ in Euclidean 4--space, the left normal $N: M \to S^2$ of
  $f$ is harmonic, hence we can associate a second Riemann surface:
  the eigenline spectral curve of $N$, as defined in
  \cite{hitchin}. We show that the multiplier spectral curve of a
  Hamiltonian stationary torus and the eigenline spectral curve of its
  left normal are biholomorphic Riemann surfaces of genus
  zero. Moreover, we prove that all Darboux transforms, which arise
  from generic points on the spectral curve, are Hamiltonian
  stationary whereas we also provide examples of Darboux transforms
  which are not even Lagrangian.
\end{abstract}

%%% Local Variables: 
%%% mode: latex
%%% TeX-master: "doc"
%%% End: 

%% file: intro.tex
\section{Introduction}

A submanifold $X$ of a real $2n$--dimensional symplectic manifold
$(Y^{2n},\omega)$ is defined to be Lagrangian if it is isotropic with
respect to the symplectic non degenerate form $\omega$ and its
dimension is maximal, that is $\dim X=n$. If furthermore $Y$ is
equipped with a Riemannian metric (in particular if $Y$ is
K\"{a}hler), a quite natural question is to ask what are the
area-minimizing Lagrangian submanifolds. In 1993 Oh \cite{oh} (after
previous work of Chen--Morvan \cite{chen_morvan}) introduced a further
variational problem intertwining Riemannian and symplectic geometry: A
Hamiltonian variation is given by a compactly supported vector field
$Z$ such that $Z \lrcorner \, \omega$ is an exact 1--form, and a
Lagrangian submanifold is called \emph{Hamiltonian stationary} if it
is a critical point of the area functional with respect to all
Hamiltonian variations. If $Y$ is K\"{a}hler--Einstein, a $S^1$-valued
function $e^{i \beta}$ can be defined along the submanifold $X$, and
$\beta$ is called the \emph{Lagrangian angle} of $X$.  It turns out
that a submanifold $X$ is Hamiltonian stationary if and only if
$\beta$ is harmonic.

In this paper we investigate immersions $f: M\to\C^2$ of a Riemann
surface $M$ into complex 2-space $\C^2$ whose images $X=f(M)$ are
Hamiltonian stationary surfaces in $\C^2$.  Even in this simplest non
trivial case Hamiltonian stationary surfaces display a rich
geometry. While no sphere can exist because of the harmonicity of
$\beta$, there exist tori, e.g. the Clifford torus $S^1 \times S^1$,
and more non--trivial examples were introduced in
\cite{castro_urbano}, \cite{anciaux}, \cite{chen}, \cite{castro_chen},
\cite{pascal_frederic}.  In \cite{pascal_frederic} it is also shown
that all Hamiltonian stationary Lagrangian tori $f: T^2\to\C^2$ are
given by Fourier polynomials, and are in this sense completely
described. On the other hand, since a Hamiltonian stationary torus $f:
T^2\to\C^2$ is given by a harmonicity condition we expect to see a
spectral curve description of $f$ when introducing a spectral
parameter.  In fact, it is possible to describe Hamiltonian stationary
tori in terms of a spectral curve \cite{pascal_frederic},
\cite{pascal_quaternionic}, \cite{pascal_ian} but the construction
shows a surprising singular behavior.

There are further notions of a spectral curve associated with a
conformally immersed Hamiltonian stationary torus: a conformal
immersion $f: M \to \R^4$ of a Riemann surface $M$ into Euclidean
4--space induces a complex structure $J$ on the trivial $\H$ bundle
$\trivial{}=\ M\times\H$ which is harmonic in the case of a
Hamiltonian stationary immersion.  In particular, one can define
\cite{uhlenbeck}, \cite{hitchin}, \cite{Klassiker}, a $\C_*$--family
of complex flat connections $d^\mu$ on the trivial $\C^2$ bundle over
$M$.  In the case when $M= T^2$ is a 2--torus, the set $\Eig$ of
eigenvalues of the holonomy of $d^\mu$ is analytic, and the
\emph{(eigenline) spectral curve} $\Sigma_e$ of the harmonic complex
structure $J$ is defined as the normalization of $\Eig$.

On the other hand, for every conformal immersion $f: T^2 \to S^4$ of a
2--torus $T^2=\C/\Gamma$ into the 4-sphere the \emph{(multiplier)
  spectral curve} is defined \cite{conformal_tori}, see also
\cite{schmidt_grinevich}, \cite{taimanov_weierstrass}: the complex
structure $J$ on $\trivial{}$ gives a quaternionic holomorphic
structure on $\trivial{}$ by $D = d''$ where $d$ is the trivial
connection on $\trivial{}$ and $d''$ denotes the $(0,1)$-part of $d$
with respect to the complex structure $J$.  We denote by $\ttrivial$
the pullback of the trivial bundle under the projection $\C
\to\C/\Gamma$, and by $ H^0_h(\ttrivial)$ the set of holomorphic
sections $\varphi\in\Gamma(\ttrivial)$ with multiplier
$h\in\Hom(\Gamma, \C_*)$, that is $D\varphi=0$ and $\gamma^*\varphi =
\varphi h(\gamma)$ for all $\gamma\in\Gamma$.  The multiplier spectral
curve $\Sigma$ of a conformal torus $f$ is defined as the
normalization of the set of all multipliers of holomorphic
sections. For a generic multiplier $h\in\Sigma$ the space of
holomorphic sections $H^0_h(\ttrivial)$ with multiplier $h$ is
1--dimensional. The lines $\L_h = H^0_h(\ttrivial)$ extend smoothly to
a line bundle, the \emph{kernel bundle}, over $\Sigma$.  Moreover, the
\emph{prolongation} of a holomorphic section $\varphi\in
H^0_h(\ttrivial)$ with multiplier $h$ defines a conformal torus $\hat
f: T^2\to S^4$ which is geometrically a \emph{Darboux transform} of
$f$. Vice versa, every closed Darboux transform gives a holomorphic
section with multiplier. In particular, the spectral curve $\Sigma$ of
a conformal torus $f$ is essentially given by the set of all closed
Darboux transforms of $f$.
 
The purpose of this paper is to study the geometry of Darboux
transforms of a Hamiltonian stationary torus $f: T^2\to\R^4$ and the
relationship between the multiplier and eigenline spectral curve of
$f$: since in this case the complex structure $J$, which is induced by
$f$ on $\trivial{}$, is harmonic and takes values in a unit circle, it
is possible to describe the set of multipliers explicitly in terms of
the lattice $\Gamma$ with $T^2=\C/\Gamma$ and the Lagrangian angle
$\beta$ of $f$. From this, we get an explicit description of all
holomorphic sections with multiplier.  In particular, for a generic
multiplier $h$ every holomorphic section $\alpha\in H^0_h(\ttrivial)$
is given by a Fourier monomial. Furthermore, this description also
yields a conceptual proof of the result of \cite{pascal_frederic} that
every lattice $\Gamma$ and $\beta_0\in\Gamma^*$ uniquely prescribes a
family of Hamiltonian stationary tori with Lagrangian angle
$\beta=2\pi\langle\beta_0,.\rangle$ provided that $\beta_0$ satisfies
a non--degeneracy condition; this follows since all Hamiltonian
stationary tori $\tilde f: T^2\to\R^4$ with Lagrangian angle $\beta$
are holomorphic sections with trivial multiplier. The family of
Hamiltonian stationary tori with the same Lagrangian angle is
therefore obtained by projection of a holomorphic curve in $\HP^k$ to
$\HP^1$ where $k=\dim_\H H^0(\trivial{})$ is the dimension of the
space of holomorphic sections.

We call holomorphic sections which are given by a Fourier monomial,
and the associated Darboux transforms, \emph{monochromatic}. We show
that all monochromatic Darboux transforms of a Hamiltonian stationary
torus are again Hamiltonian stationary, and that in this case the
Lagrangian angle is, up to reparametrization, preserved.  The space of
holomorphic sections with a given multiplier is generically complex
1--dimensional, and only at multiplier with high--dimensional space of
holomorphic sections polychromatic holomorphic sections may occur. In
particular, for all but a finite set of multipliers, we see that the
associated Darboux transforms are Hamiltonian stationary.  Discussing
the example of homogeneous tori and Castro Urbano tori, we show that
there exist however families of non--Lagrangian, polychromatic Darboux
transforms, each family associated with a multiplier with
high--dimensional space of holomorphic sections.

On the other hand, the family of flat connections $d^\mu$ on the
trivial $\C^2$--bundle gives rise to a subset of Darboux transforms,
the so--called $\mu$--Darboux transforms, since it turns out that
every parallel section of $d^\mu$ is indeed holomorphic. We show that
all (even local) $\mu$--Darboux transforms of a Hamiltonian stationary
immersion $f: M \to \R^4$ have harmonic left normal and are thus
constrained Willmore.  Furthermore, we see that in the case of a
Hamiltonian stationary torus $M =T^2$, all $\mu$--Darboux transforms
$\hat f: T^2\to\R^4$ are given by monochromatic holomorphic sections,
and thus are Hamiltonian stationary. Conversely, every monochromatic
holomorphic section $\alpha$ with multiplier gives rise to a unique
$\mu\in\C_*$ such that $d^\mu\alpha=0$. In particular, this
correspondence induces a biholomorphism from the eigenline spectral
curve to the multiplier spectral curve.  This way, we see that the
multiplier spectral curve can be compactified to a connected Riemann
surface of genus zero. Though the normalizations of the multiplier
spectral curve and the spectral curve in \cite{pascal_ian} are 
related this also shows that they do not coincide. However, the
original harmonic complex structure and the conformal torus can
already be recovered from the spectral curve $\Sigma$ and its kernel
bundle by appropriate limits as $\mu\to\infty$.

%%% Local Variables: 
%%% mode: latex
%%% TeX-master: "doc"
%%% End: 

%% file: holomorphic.tex
\section{Holomorphic sections with multiplier}
\label{sec: holo sections w. multiplier}

To discuss the multiplier spectral curve of a Hamiltonian stationary
torus we briefly recall the general construction for conformal tori
\cite{conformal_tori}: In the following we will identify a conformal
immersion $f\colon M \to S^4$ with the quaternionic line subbundle $L
\subset V$ of the trivial $\H^2$ bundle $V = \trivial{2}$ whose fibers
are $L_p = f(p)\in \HP^1$ for $p\in M$.  All quaternionic vector spaces
are here, and in what follows, quaternionic right vector spaces. Since
$f$ is a immersion, the derivative $\delta=\pi d|_L$ of $f$ has no
zeros where $d$ is the trivial connection on $V$, and $\pi\colon V \to
V/L$ is the canonical projection. In particular, the conformality of
$f$ defines \cite[Section 2.5]{Klassiker}  a complex structure
$J\in\Gamma(\End(V/L))$, $J^2=-1$, on $V/L$ so that
\begin{equation}
\label{eq:conformality}
*\delta = J \delta\,.
\end{equation}
The complex structure $J$ and the trivial connection $d$ on $V$ induce
\cite{conformal_tori} a quaternionic holomorphic structure $D\colon
\Gamma(V/L) \to\Gamma(\bar K V/L)$ on $V/L$ via
\[
D\varphi = (\pi d\tilde\varphi)''
\]
where $\tilde \varphi\in\Gamma(V)$ is an arbitrary lift of
$\varphi\in\Gamma(V/L)$, that is $\pi\tilde\varphi =
\varphi$. Moreover, we denote by
\[
\omega' = \frac{1}{2}(\omega - J *\omega) \quad \text{ and } \quad \omega'' = \frac{1}{2}(\omega + J *\omega)
\]
the $(1,0)$ and $(0,1)$--part of a 1--form $\omega\in\Omega^1(E)$ with
values in a complex vector bundle $(E, J)$.  Note that $D$ is
independent of the choice of $\tilde \varphi$ since
$\pi(d\tilde\psi)'' =(\delta\tilde\psi)''=0$ for all
$\tilde\psi\in\Gamma(L)$ by (\ref{eq:conformality}).  We call
$\varphi\in \ker D$ a \emph{(quaternionic) holomorphic section}, and
denote by $H^0(V/L) = \ker D$ the space of holomorphic
sections. Furthermore, we write $\widetilde{V/L}$ for the pullback of
$V/L$ under the projection $\tilde M \to M$ of the universal cover
$\tilde M$ to $M$.
\begin{definition}[see \cite{conformal_tori}]
  If $\varphi\in H^0(\widetilde{V/L})$ with $\gamma^* \varphi =\varphi
  h_\gamma$ and $h_\gamma \in\H_*$ for all $\gamma\in\pi_1(M)$, then
  we call $h: \pi_1(M) \to\H_*$ the \emph{multiplier} of the
  holomorphic section $\varphi$. We denote by $ H^0_h(\widetilde{V/L})
  $ the set of holomorphic sections with multiplier $h$.
\end{definition}
If $M = T^2$ is a 2--torus with lattice $\Gamma$ the image
of $h\colon \Gamma \to \H_*$ lies in an abelian subgroup of $\H_*$.
By scaling $\varphi$ with a quaternion we may thus assume that
$h$ takes values in $\C= \Span\{1,i\}$. If $f$ has zero normal
bundle degree the set of multipliers
\[
\Spec =\{ h \colon \Gamma \to \C_*\mid \text{ there exists } \varphi\in H^0(\widetilde{V/L}) \text{ with } \gamma^*\varphi = \varphi h_\gamma \text{ for all 
} \gamma\in\Gamma
\}
\]
is a 1-dimensional analytic variety, and the normalization of $\Spec$
is a Riemann surface $\Sigma$ with at most two connected components,
each containing a point at infinity, and with possible infinite genus.
The normalization $\Sigma$ of $\Spec$ is called the \emph{multiplier
  spectral curve} of the conformal torus $f\colon T^2\to S^4$.  Since
\[
\gamma^*(\varphi j)
= (\varphi j)\bar h_\gamma
\]
for a holomorphic section $\varphi\in H^0(\widetilde{V/L})$ with
multiplier $h$, we see that $h\in\Spec$ implies $\rho(h)\assign\bar
h\in\Spec$. In fact, $\rho$ is induced by a fixed point free real
structure $\rho$ on the multiplier spectral curve $\Sigma$.

Moreover, there is a complex holomorphic line bundle $\L$ defined over
the multiplier spectral curve, the so--called \emph{kernel bundle}: at
a generic point $h\in\Sigma$ the space of holomorphic sections with
multiplier $h$ is 1--dimensional, and defines a complex line
\[
\L_h:=H^0_h(\widetilde{V/L})\subset\Gamma(\widetilde{V/L})\,.
\] 
The lines $\L_h$ extend smoothly into points on the multiplier
spectral curve with high--dimensional space of holomorphic sections,
and thus define a holomorphic line subbundle $\L$ of the trivial
$\Gamma(\widetilde{V/L})$ bundle over $\Sigma$.  If $\Sigma$ has
finite genus, then $\Sigma$ can be compactified and $\L$ extends
smoothly to the compactified spectral curve $\bar \Sigma$. Moreover,
the line bundle $\L$ is compatible with the real structure $\rho$ that
is $\rho^*\L = \L j$.  

We now turn to the case when $f\colon T^2\to \R^4$ is a conformal immersion
into the 4--space with Gauss map
\[
(N, R)\colon M \to S^2 \times S^2 = \Gr_2(\R^4)\,.
\]
If we identify the Euclidean 4--space $\R^4=\H$ with the quaternions
then $S^2 = \{ n\in\Im\H \mid n^2=-1\} $ and the \emph{left normal}
$N$ and the \emph{right normal} $R$ satisfy
\[
*df = N df = -df R\,.
\]
We consider a conformal immersion $f\colon M \to\R^4$ as a map into
$S^4=\HP^1$ via $\psi\H\colon M \to\HP^1$ where
\[
\psi=\begin{pmatrix} f\\1
\end{pmatrix}\colon M \to \H^2\,.
\]
In other words, the map $\psi\H\colon M \to S^4$ becomes $f\colon M \to\R^4$
after the choice of the point $\infty=e\H \in S^4$ where $e=\begin{pmatrix} 1\\0
\end{pmatrix}$.  In the case of a conformal immersion $f\colon M
\to\R^4$ by evaluating $\delta $ on $\psi$ we see that the complex
structure $J$ on $V/L$ is given by $J\pi e= \pi e N$ where $N$ is the
left normal $N$ of $f$.  In what follows we identify
\[
V/L \cong \underline{e\H}
\]
via the quaternionic isomorphism $\pi|_{\underline{e\H}}\colon
\underline{e\H} \to V/L, e \mapsto \pi(e)$, and trivialize
$\underline{e\H} \cong \trivial{}$ via the constant section $e$. In
particular, $\trivial{}$ inherits the complex structure
$J\in\End(\trivial{})$ which is given by left multiplication by $N$.
Moreover, the holomorphic structure $D$ on $\trivial{}$ is given by
\[
D\alpha = \frac{1}{2}(d\alpha + N*d\alpha)
\]
so that $\alpha\in H^0(\trivial{})$ is holomorphic if and only if 
\begin{equation}
\label{eq:holomorphicity}
*d\alpha = N d\alpha\,.
\end{equation}

We shall consider Lagrangian surfaces of $\R^4=\H$ where the complex
structure on $\H$, which determines the symplectic structure, is given
by the left multiplication by $j$. In other words, we identify $\H =
\Cj \oplus \Cj i$ where $\Cj =\Span\{1,j\}$ so that the action of a
unitary matrix
\[
U = e^{j\theta} \begin{pmatrix} m & -\bar n\\ n &\bar m
\end{pmatrix}
\in U(2), \quad m, n\in \Cj, \theta\in\R\,,
\] 
on $\Cj^2$ corresponds to left multiplication with $e^{j\theta}$
together with right multiplication by $m+ n i$ on $\H$:
\[ 
\Cj^2 =\H \ni v \mapsto U v = e^{j\theta} v (m + ni)\in \H\,.
\]
A conformal immersion 
 $f:M\to\H$ is Lagrangian
 if and only if
\[
\left(e^{-u} \frac{\partial f}{\partial x}, e^{-u} \frac{\partial f}{ 
\partial y}\right)
\]
is a unitary frame where $z = x + iy$ is a local conformal coordinate
on $M$ and $e^{u}$ is the conformal factor of $f$. Choosing the
couple $(1,i)$ as reference unitary frame, there exists locally a
smooth $U \colon M \to U(2)$ such that $\frac{\partial f}{\partial x}=U\cdot 1$ and
$\frac{\partial f}{\partial y}= U\cdot i$, that is there exist $\beta \colon M\to
\R$ and $q \colon M\to S^3$ such that
\[
df = \frac{\partial f}{\partial x} dx+\frac{\partial f}{\partial y}  
dy = e^{\frac{j \beta}2} (1 dx+i dy ) e^u q\,.
\]
Writing $g = e^u q$ we thus see that 
  $f\colon M\to\H$ is conformal and Lagrangian if and only if
\begin{equation}
\label{eq:df}
df = e^{ \frac{j\beta} 2} dz g\; , \textrm{ where } dz=dx+i dy,
\end{equation}
for some real valued function $\beta$ and quaternion valued $g$.  If
we change the conformal coordinate $\tilde z = z e^{i\theta}$ we see
that $ \tilde g = e^{-i\theta} g$. In particular, the \emph{Lagrangian
  angle} $\beta$ is defined independently, up to a constant
translation, of the choice of the conformal coordinate $z$. It is now
easily seen from (\ref{eq:df}) that the left normal is
\begin{equation}
\label{eq:N}
N=e^{j \beta} i
\end{equation}
while the right normal is 
\begin{equation}
\label{eq:right normal}
R=  -g^{-1} i g\,,
\end{equation}
where we use the convention that $*dz=*(dx+i dy)=i dz = -dy + i
dx$. Conversely, if $f: M \to\R^4$ is a conformal (branched) immersion
with left normal $N = e^{j\beta} i$ with $\beta: M \to\R$ then $df$
satisfies (\ref{eq:df}) for some $g: M \to\H$, and thus $f$ is a
Lagrangian immersion.  Since our interest are Hamiltonian stationary
immersions from a 2--torus $T^2=\C/\Gamma$ into $\R^4$, we will
consider $f\colon \C \to \H$ as a $\Gamma$-periodic map. In
particular, if $f$ is Hamiltonian-stationary, then $\beta$ is harmonic
and we may assume, after possible change of the conformal coordinate
$z$, that
\begin{equation}
\label{eq:beta_0}\beta (z) =
2 \pi \left\langle \beta_0, z \right\rangle\,,
\end{equation}
where $\left\langle, \right\rangle$ is the scalar product in $\C
\simeq \R^2$ and $\beta_0$ is in the dual lattice $\Gamma^{\ast}$.
 
In the following we discuss the set of all multipliers of holomorphic
sections in case of a Hamiltonian stationary torus.  Recall that the
multiplier takes values in any complex subspace of $\H$ since the
image of $h: \Gamma \to \H_*$ is abelian. We choose this subspace to
be $\C =\Span\{1, i\}$ instead of $\Cj =\Span\{1, j\}$ for purely
computational reasons. In particular, for every $h\in\Spec$ there exists
a pair $(A,B) \in\C^2$ so that
\[
h_\gamma =h^{A,B}_\gamma \assign e^{2\pi(\langle A,\gamma\rangle -i\langle B, \gamma\rangle)}\, \quad \gamma\in\Gamma\,.
\]
Note that $(A,B)$ is unique  up
to translation of $B$ by elements of the dual lattice since
\begin{equation}
\label{eq:lattice invariance}
h^{A,B}  = h^{A, B+\delta} \quad \text{ for all } \delta \in\Gamma^*\,.
\end{equation}
To find conditions on $(A,B)$ for $h^{A,B}$ to be the multiplier of a
holomorphic section $\alpha$ we rewrite the holomorphicity condition
in terms of the gauged section $\tilde{\alpha} \assign e^{-
  \frac{j\beta} 2} \alpha$. In particular, using $\frac{1}2(d \beta +
i *d \beta) = \pi \beta_0 d \bar{z}$ and (\ref{eq:N}),
\begin{equation}
  d\tilde{\alpha} + i \ast d \tilde{\alpha} + \pi \beta_0 d\bar z j \tilde{\alpha} = 0\,,
  \label{gaugedholomorphicity}
\end{equation}
 is equivalent to (\ref{eq:holomorphicity}).
  If $\alpha\in H^0_h(\widetilde{V/L})$ is a
holomorphic section with multiplier $h=h^{A,B}$ then $\alpha (z) =
\sigma e^{2 \pi ( \left\langle A, z \right\rangle - i \left\langle B,
    z \right\rangle)}$ with  $\Gamma$-periodic function  $\sigma$ so that 
\begin{equation}
\label{eq:tilde alpha}
\tilde \alpha (z) = \tilde \sigma e^{2 \pi ( \left\langle A, z
  \right\rangle - i \left\langle B, z \right\rangle)}
\end{equation}
where $\tilde \sigma =\emjbh \sigma$ is $2\Gamma$--periodic. Writing
$\tilde \sigma = \tilde u + j\tilde v$ with complex valued functions
$\tilde u, \tilde v$ equation (\ref{gaugedholomorphicity}) then
becomes
\begin{eqnarray*}
  0 & = &  \{2 \tilde u_{\bar{z}} + 2 \pi (A - iB) \tilde u - \pi \beta_0 \tilde v\} d\bar{z} + j
  \left\{ 2 \tilde v_z + 2 \pi ( \bar{A} - i \bar{B}) \tilde v + \pi \bar{\beta}_0 \tilde u
  \right\} d z\,,
\end{eqnarray*}
and we have to discuss for which $(A,B)\in\C^2$ there exist complex valued $2\Gamma$--periodic
functions  $\tilde u, \tilde v$ satisfying 
\begin{equation}
  \left\{ \begin{array}{l}
    2 \tilde u_{\bar{z}} + 2 \pi (A - iB) \tilde u - \pi \beta_0 \tilde v = 0\\[.2cm]
    2 \tilde v_z + 2 \pi ( \bar{A} - i \bar{B}) \tilde v + \pi \bar{\beta}_0\tilde u = 0
  \end{array}\,. \right. \label{holouv}
\end{equation}
Since $\sigma = u + j v\colon T^2\to\H$ is defined on the torus  we have  Fourier expansions
\[
u(z) = \sum_{\delta\in\Gamma^*} u_\delta \edelta \quad \text{ and  } \quad
v(z) = \sum_{\delta\in\Gamma^*} v_\delta \edelta
\]
with $u_\delta, v_\delta\in\C$. Abbreviating $e_\delta(z) = \edelta$ we obtain with (\ref{eq:beta_0})
\[
\emjbh = \cos\frac\beta 2  - j\sin \frac{\beta}{2} 
= \frac{\ebnh + \embnh}2 - i j \frac{\ebnh - \embnh}{2}
 \]
so that $\tilde u+j\tilde v = \emjbh(u +
jv )$ and $j \ebnh = \embnh j$ give 
  the Fourier expansions
\begin{equation}
\label{eq:Fourier expansion}
\tilde u = \sum_{\delta\in \Gamma^*+\frac{\beta_0}2}\tilde u_{\delta}
e_\delta\quad  \text{ and  } \quad 
\tilde v = \sum_{\delta\in \Gamma^*+\frac{\beta_0}2}\tilde  v_{\delta} e_\delta 
\end{equation}
of $\tilde u$ and $\tilde v$ with Fourier coefficients $\tilde u_\delta, \tilde
v_\delta\in\C$. Therefore, the gauged holomorphicity equations
(\ref{holouv}) can be written equivalently as
\begin{equation}
  \left\{ \begin{array}{l}
    2 (i \delta + A - iB) \tilde u_{\delta} = \beta_0  \tilde v_{\delta} \\[.2cm]
    2 (i \bar{\delta} + \bar{A} - i \bar{B}) \tilde v_{\delta} =-
    \bar{\beta}_0 \tilde  u_{\delta} 
  \end{array} \right. \label{Fourieruv} \quad \text{ for all
 }\delta\in\Gamma^*+\frac{\beta_0}2\,.
\end{equation}
In particular, the vanishing of one of the Fourier coefficients
implies the vanishing of the other since $\beta_0 \neq 0$.  For
$\tilde u_{\delta} \tilde v_{\delta} \neq 0$ the above equations imply
\begin{equation*}
(\delta - iA - B) ( \bar{\delta} - i \bar{A} - \bar{B}) = \frac{| \beta_0
   |^2}{4}\,, 
\end{equation*}
or, equivalently, 
\begin{equation}
  | \delta - B|^2 - |A|^2 = \frac{| \beta_0 |^2}{4}   
\quad \text{ and } \quad \langle \delta-B, A\rangle =0.
\label{ABcondition}
\end{equation}
In particular, the condition that the set of admissible frequencies
 \[ \Gamma^*_{A,B}
=\{\delta\in\Gamma^*+\frac{\beta_0}2 \mid \delta \text{ satisfies
  (\ref{ABcondition})}\}
\]
is not empty is necessary for $h^{A,B}$ to be a multiplier of a
holomorphic section. \\
  
\begin{theorem}
 \label{thm:parametrized multipliers}
 A multiplier $h\in\Spec$ of a holomorphic section $\varphi\in
 H^0(\widetilde{\trivial{}})$ determines a unique pair $(A,B)\in\C^2$,
 up to translation of $B$ by the dual lattice $\Gamma^*$, with
 $h^{A,B} = h$.  Moreover, the set of multipliers of a Hamiltonian
 stationary torus $f: T^2 \to S^4$   is given by
\[
\Spec = \{ h^{A,B} \mid  \Gamma^*_{A,B}\not=\varnothing
\}\,.
\]
In fact, for all $\delta\in\Gamma^*_{A,B}$ there exists a complex
1--dimensional subspace $\L_{A,\delta-B}$ of the space
$H^0_{A,B}\assign H^0_{h^{A,B}}(\widetilde{\trivial{}})$ of
holomorphic sections with multiplier $h^{A,B}$.
\end{theorem}

\begin{proof} 
  We already have seen that 
  $\Spec\subset \{ h^{A,B} \mid \Gamma^*_{A,B}\not=\varnothing\}$.
  To show equality, note that the
  frequencies $\delta$ with $\tilde u_\delta \tilde v_\delta \not=0$
  are placed at the intersection of the circle $\mathcal{C}_r(B)$
  of radius 
\[
r =\sqrt{|A|^2 + \frac{|\beta_0 |^2}4}
\]
centered at
  $B$, and a line going through $B$ orthogonal to $A$. If $A = 0$, the
  second condition is void and $(0, B)$ satisfies (\ref{ABcondition})
  if and only if the circle $\mathcal{C}_{\frac{|\beta_0|}2}(B)$ meets
  $\Gamma^{\ast} + \frac{\beta_0}{2}$. In other words,
 if $\delta\in\Gamma^*_{0,B}$ then $\delta = B -
  \frac{\beta_0}{2} e^{it}$ for some $t\in[0,2\pi)$, that is
\begin{equation}
\label{eq: Gamma 0B}
\Gamma^*_{0,B} = \{\delta \in \Gamma^* + \frac{\beta_0}2 \mid \delta= B -\frac{\beta_0}2e^{it}, t\in [0,2\pi)\}\,.
\end{equation}
In particular, for $\delta\in\Gamma^*_{0,B}$ we have
\begin{equation*}
\frac{i\bar\beta_0}{2(\bar\delta -\bar B) }= \frac{2i(\delta-B)}{\beta_0} = -ie^{it}\,,
\end{equation*}
so that  (\ref{eq:tilde alpha}), (\ref{eq:Fourier expansion}), and
~(\ref{Fourieruv}) show that
\begin{equation}
\label{eq:alpha_gamma, A=0}
\alpha_{\delta} = e^{\frac{j\beta}2}  (1 + k e^{it})e_{-\frac{\beta_0 e^{it}}2}
\end{equation}
is a holomorphic section with multiplier $h^{0, B}$. In particular,
$\alpha_\delta$ spans a complex 1--dimensional subspace
$\L_{A,\delta-B}\subset H^0_{0,B}$ of the space of holomorphic
sections with multiplier $h^{0,B}$.

If $\delta\in \Gamma^*_{A,B}$ with $A\not=0$, then $\delta$
lies on the intersection of a circle of radius $ r =\sqrt{|A|^2 +\frac{ | \beta_0 |^2}4}$ centered at
$B$ and a line going through $B$ orthogonal to $A$, that is 
\begin{equation}
\label{eq:Gamma_A,B, else}
\Gamma^*_{A,B} =  \left\{
\delta_\pm=  B  \pm  i r\frac{A}{|A|} \right\} \cap \left(\Gamma^* + \frac{\beta_0}2\right)\,.
\end{equation}
By (\ref{ABcondition}) we see that 
\begin{equation*}
\frac{-\bar\beta_0}{2(i\bar\delta_\pm + \bar A - i\bar B)} = \frac{2(i\delta_\pm + A - iB)}{\beta_0} 
\end{equation*}
and, if $\delta_\pm \in\Gamma^*_{A,B}$, again  (\ref{eq:tilde alpha}), (\ref{eq:Fourier expansion}), and
~(\ref{Fourieruv})  together with
\[
 \frac{2(i\delta_\pm + A - iB)}{\beta_0} = \frac{2(|A| \mp r)}{\beta_0}\frac A{|A|}
\]
yield that
\begin{equation}
\label{eq: alpha_gamma, else}
\alpha_{\delta_\pm} =  e^{\frac{j\beta}2}  \left(1 + j \frac{2( |A| \mp r)}{\beta_0} \frac A{|A|} \right) e_{\delta_\pm} e^{2\pi(\langle A, \cdot \rangle - i\langle B, \cdot \rangle)} 
\end{equation}
is a holomorphic section with multiplier $h^{A,B}$, and  $\L_{A,\delta_\pm -B}\subset H^0_{A,B}$ is spanned by $\alpha_{\delta_\pm}$.
\end{proof}

We call a holomorphic section $\alpha\in H^0_{A,B}$
\emph{monochromatic} (respectively \emph{polychromatic}) if $\alpha$
is given, up to the gauge, by a Fourier monomial (respectively by a
Fourier polynomial with more than one frequency). As we have seen, any
monochromatic holomorphic section is given by (a complex scale of)
\begin{equation}
\label{eq:monochromatic holomorphic section}
\alpha_\delta=\ejbh (1-k\lambda_\delta)  e_{\delta-B} e^{2\pi\langle A, \cdot\rangle}
\end{equation}
with $\delta\in\Gamma^*_{A,B}$  and (\ref{eq:alpha_gamma, A=0}), (\ref{eq: alpha_gamma, else})
\begin{equation}
\label{eq:lambda delta}
\lambda_\delta \assign \lambda_{A, \delta- B}\assign \frac{2}{\beta_0}(\delta-i A - B) =\begin{cases} -e^{it} \quad & \text{for } A =0, \delta = B -\frac{\beta_0}2e^{it},\\
-\frac{2i(|A| \mp r)}{\beta_0}\frac{A}{|A|} \quad &\text{for } A\not=0, \delta =\delta_\pm\,
\end{cases}
\end{equation}
where $r=\sqrt{|A|^2 +\frac{|\beta_0|^2}{4}}$.  Since every multiplier
$h$ is given by $h=h^{A,B}$ for some pair $(A,B)$ satisfying
(\ref{ABcondition}), we also obtain:
\begin{cor}
\label{cor:holomorphic sections}
  Every holomorphic section $\alpha\in H^0_h(\widetilde{\trivial{}})$
  with multiplier $h$ is given by
\[
\alpha = \ejbh\left(\sum_{\delta\in\Gamma^*_{A,B}}(1-k\lambda_\delta)\tilde u_\delta e_{\delta-B}\right) e^{2\pi\langle A, \cdot\rangle}\,,
\]
where $h=h^{A,B}$, $\tilde u_\delta\in\C$, $e_{\delta-B}(z)=\edeltaB$,  and 
$
\lambda_\delta= \frac{2}{\beta_0}(\delta-iA -B)\,.
$ 
\end{cor}\begin{rem}
  Note that $\alpha$ is independent of the choice of $B$ with
  $h^{A,B}=h$ since both $\lambda_\delta$ and $e_{\delta-B}$ only
  depend on $\delta-B$, and $\Gamma^*_{A,B + \zeta} =\Gamma^*_{A,B}
  +\zeta$ for all $\zeta\in\Gamma^*$. 
\end{rem}

Since $H^0_{A,B}$ is complex $k$--dimensional if and only if the
Fourier expansion (\ref{eq:Fourier expansion}) allows exactly $k$
frequencies $\delta\in\Gamma^* + \frac{\beta_0}2$ with $\tilde
u_{\delta} \tilde v_{\delta}\not=0$,   the
complex dimension of the space of holomorphic sections with multiplier
$h^{A,B}$ is given by the number of elements in $\Gamma^*_{A,B}$:
\[
\dim_\C H^0_{A,B} = |\Gamma^*_{A,B}|\,.
\]
Note that $A\not=0$ uniquely determines $B$, up to translation by the
dual lattice, with $\Gamma^*_{A,B}\not=\varnothing$, and
(\ref{eq:Gamma_A,B, else}) shows that $|\Gamma^*_{A,B}|\le 2$ when
$A\not=0$. If $|\Gamma^*_{A,B}| = 2$ then both $\delta_\pm = B \pm i
r\frac{A}{|A|} \in\Gamma^* + \frac{\beta_0}2$ are in the translated
dual lattice (\ref{eq:Gamma_A,B, else}) so that
\[
 \delta_+ - \delta_- =: \zeta \in\Gamma^*\,.
\]
On the other hand, $ \zeta = 2iA \sqrt{1 + \frac{|\beta_0|^2}{4|A|^2}}
$ gives $4|A|^2 = |\zeta|^2 -|\beta_0|^2$, and $A = -\frac{i}2
\zeta \sqrt{1 - \frac{|\beta_0|^2}{|\zeta|^2}}$ with $|\zeta| >
|\beta_0|$. \\

\begin{lemma}
 Let $f$ be a Hamiltonian stationary torus and   $h=h^{A,B}\in\Spec$ a multiplier of a holomorphic section with $A\not=0$.
Then the complex dimension of the space of holomorphic sections with multiplier $h$ is at most 2, and $\dim H^0_{A,B}=2$ if and only if $A$ is a \emph{double point} that is
\[
A =  -\frac{i}2 \zeta \sqrt{1 - \frac{|\beta_0|^2}{|\zeta|^2}}, \quad \zeta\in\Gamma^*, \quad |\zeta| >
|\beta_0|\,.
  \]
\end{lemma}

For $A=0$ the situation is more complicated but we still obtain an
upper bound for the dimension of the space of holomorphic sections
with multiplier.

\begin{theorem}
  \label{thm/lemm:holomorphic sections} For a Hamiltonian stationary torus,
  the space $H^0_h(\widetilde{\trivial{}})$ of holomorphic sections
  with multiplier $h$ is generically complex 1--dimensional, and its complex dimension is bounded by $N$ where 
\[
N = |\mathcal D \cap (\Gamma^* + \frac{\beta_0}2)|
\]
is the number of points of the translated dual lattice $\Gamma^* +
\frac{\beta_0}2$ in the (closed) disk $\mathcal D$ around
$\frac{\beta_0}2$ of radius $|\beta_0|$. Moreover, there exists
at least one multiplier $h$, namely the trivial multiplier, with
$\dim_\C H^0_h(\ttrivial)\ge 4$.
\begin{figure}[h]
\label{fig:1}
  \includegraphics[height=90mm]{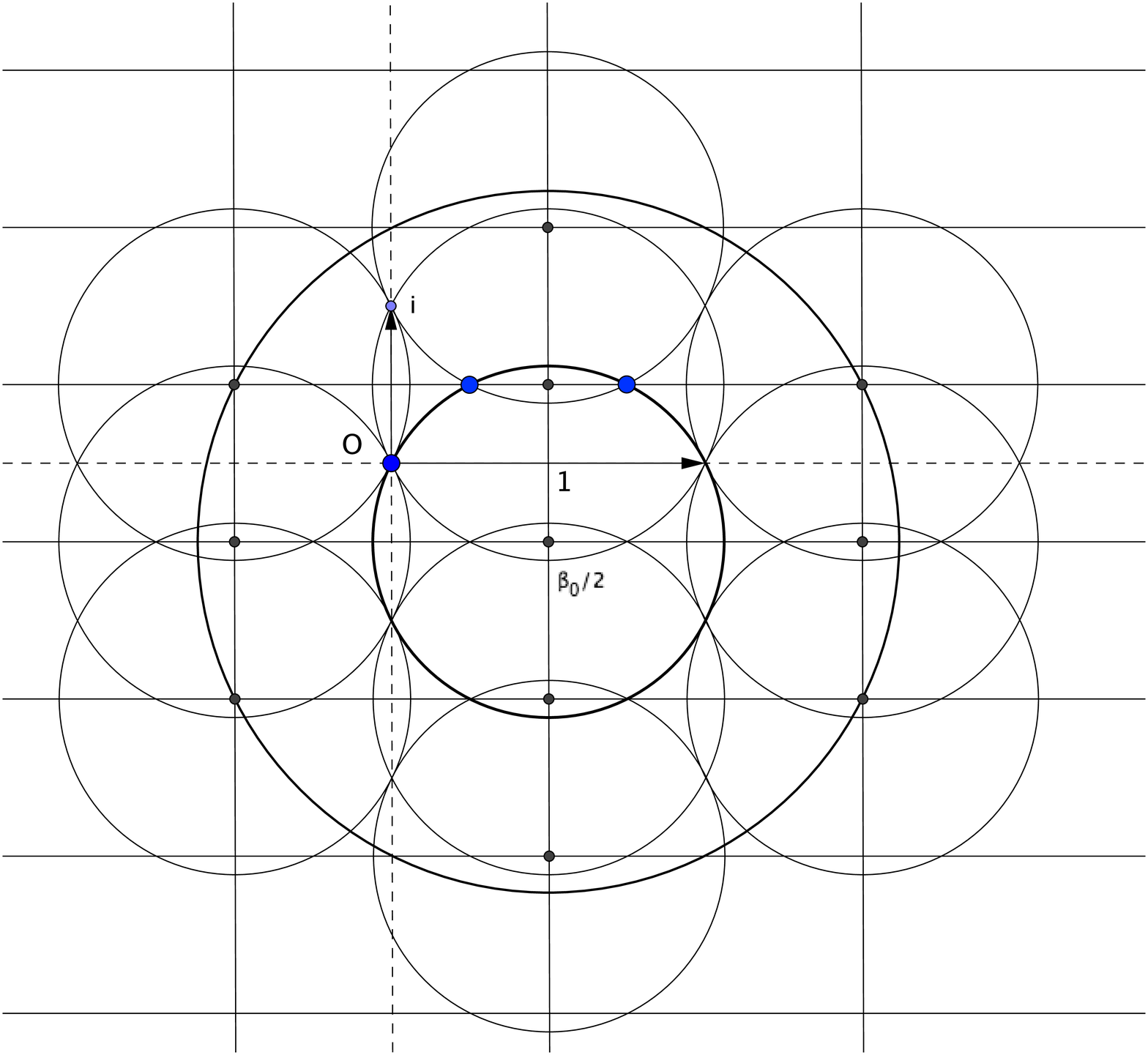}
  \caption{The points $B$ on the inner circle of radius
    $\frac{|\beta_0|}2$ around $\frac{\beta_0}2$ give all multipliers
    $h=h^{0,B}$. Frequencies $\delta\in\Gamma^*+\frac{\beta_0}2$ in
    the disk of radius $|\beta_0|$ are all frequencies such that the
    circle $\mathcal{C}_{\frac{|\beta_0|}2}(\delta)$ meets the inner
    circle.  The corresponding dimension $k$ of the space of
    holomorphic section with multiplier $h$ is $k=2$ or $k=4$, when 2
    respectively 4 circles meet at $B$, and $k=1$ for the remaining
    points. All points $B$ with $k>1$ are congruent modulo the dual
    lattice to one of the three blue points. This setting corresponds
    to a homogeneous torus, see Section 4.}
\end{figure}

\end{theorem}

\begin{proof}
  In Theorem \ref{thm:parametrized multipliers} we have seen that the
  set of multipliers of holomorphic sections is parametrized by pairs
  $(A,B)$ with $\Gamma^*_{A,B}\not=\varnothing$.  For $A\not=0$ we see
  that the complex dimension of the space of holomorphic sections with
  multiplier $h^{A,B}$ is generically 1--dimensional since $\dim_\C
  H^0_{A,B}=2$, $A\not=0$, occurs only for a discrete set of double
  points in $\C_*$.  If $A=0$, then by (\ref{eq: Gamma 0B}) 
\[
\Gamma^*_{0,B} = \{\delta \in \Gamma^* + \frac{\beta_0}2 \mid \delta= B -\frac{\beta_0}2e^{it}, t\in [0,2\pi)\}\not=\varnothing\,,
\]
and we may assume without loss of generality that $B =
\frac{\beta_0}2(1 + e^{it})$ with $t\in [0,2\pi)$, and thus $
\frac{\beta_0}2\in\Gamma^*_{0,B}$, since the multiplier $h=h^{0,B}$
does not change under translation of $B$ by the dual lattice
(\ref{eq:lattice invariance}).  Now $\delta$ lies in $\Gamma^*_{0,B}$
if and only if $\delta \in \Gamma^*+\frac{\beta_0}{2}$ and the circle
$\mathcal{C}_{\frac{|\beta_0|}2}(\delta)$ of radius
$\frac{|\beta_0|}2$ centered at $\delta$ passes through $B$.
Obviously such a $\delta$ lies inside the closed disk $\mathcal{D}$
centered at $\frac{\beta_0}2$ of radius $|\beta_0|$.  Only finitely many
points $N$ of the (translated) dual lattice are contained in
$\mathcal{D}$ which gives the upper bound
\[
\dim_\C H^0_{0,B} \le N <\infty\,.
\]
Furthermore the complex dimension of $H^0_{0,B}$ is the number of such
circles passing through $B$, which is one except for finitely many
values of $B$.
 
By (\ref{eq:holomorphicity}) we see that $1\in H^0(\trivial{})$ and $f
\in H^0(\trivial{})$ are quaternionic independent holomorphic sections
with trivial multiplier. Since the space of holomorphic sections
$H^0(\trivial{})$ is quaternionic, we thus have $\dim_\C H^0_{h\equiv
  1}(\widetilde{\trivial{}}) =2\dim_\H H^0(\trivial{}) \ge 4$.
\end{proof}
\begin{rem} Note that $\delta\in\Gamma^*_{0,B}$ if and only if the
  translate of the main circle
  $\mathcal{C}_{\frac{|\beta_0|}2}(\frac{\beta_0}2)$ by the frequency
  $\delta-\frac{\beta_0}2\in\Gamma^*$ intersects the main circle at
  $B$. Put differently, the complex dimension of the space of
  holomorphic sections with multiplier $h=h^{0,B}$ is given by the
  number of self intersections on a fundamental domain of the main
  circle at $B$.
\begin{figure}[h]
  \includegraphics[height=60mm]{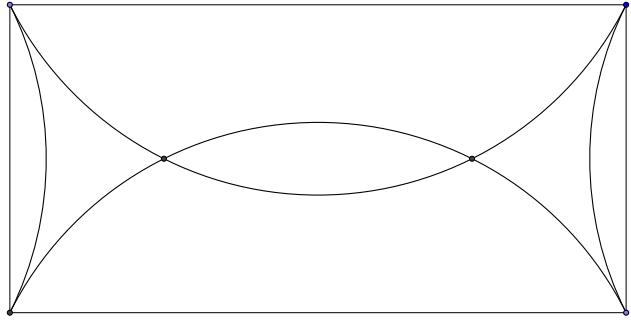}
  \caption{In the previous example, on a fundamental domain of $\C/
    \Gamma^{\ast}$ the inner circle $\mathcal{C}_{\frac{|\beta_0 | } 2}
    (\frac{\beta_0} 2)$ self intersects at three points, one of which is
    covered 4 times.}
\end{figure}
\begin{figure}[h]
  \includegraphics[height=75mm]{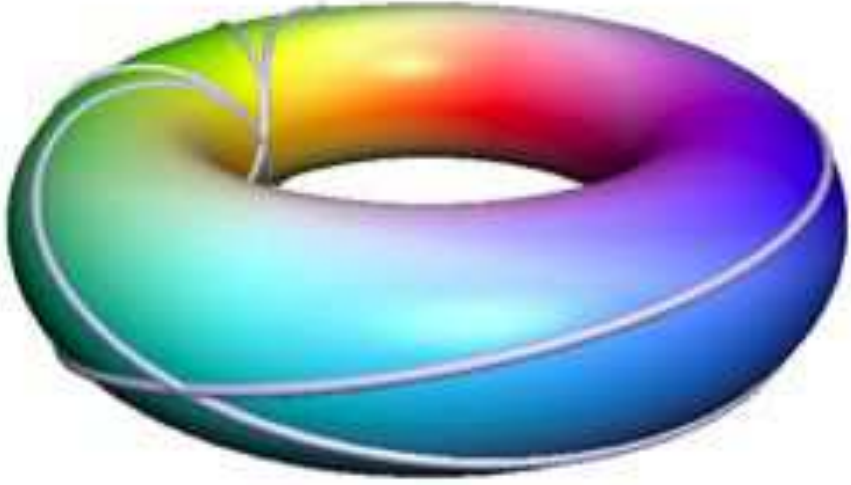}
  \caption{Same  as Figure 2 but the torus $\C/
    \Gamma^{\ast}$ is drawn in $\R^3$. Only the point with four
    intersections and one of the two double points are visible.}
\end{figure}
\end{rem}
\begin{rem}
\label{rem:trivial multiplier}
In \cite{pascal_frederic} it is shown that all Hamiltonian-stationary
tori $f: M \to\R^4$ with lattice $\Gamma$ and Lagrangian angle
frequency $\beta_0$ satisfy
\[
 \Gamma^{\ast}_{\beta_0} \assign \left\{ \delta \in \Gamma^{\ast} +
  \frac{\beta_0}{2}, | \delta | = \frac{| \beta_0 |}{2}, \delta \neq
  \pm \frac{\beta_0}{2} \right\} \not=\varnothing\,.
\]
By definition $\Gamma^*_{\beta_0}\cup \{\pm \frac{\beta_0}2\} =
\Gamma^*_{0,0}$, and $| \Gamma^{\ast}_{\beta_0}|$ is even by
symmetry. Since $h^{0,0}$ is the trivial multiplier, this also implies
that $\dim_\C H^0_{0,0} = |\Gamma^*_{0,0}| = |\Gamma^{\ast}_{\beta_0}|
+ 2 \geqslant 4$.  Since the holomorphic sections $\alpha_{\pm
  \frac{\beta_0}2} = 1 \mp k$ given by the frequencies $\pm
\frac{\beta_0}2$ are constant and $f\in H^0(\trivial{})$ the
Hamiltonian stationary torus $f$ is given up to rotation and
translation by any holomorphic section $\alpha_\delta$ with
$\delta\in\Gamma^*_{\beta_0}$ if $\dim_\C H^0(\trivial{})=4$.
\end{rem}

More generally, we can construct all Hamiltonian stationary tori
with Lagrangian angle $\beta$:

\begin{theorem}[see \cite{pascal_frederic}]
  \label{thm:hslLagrangianangle} If $\, \Gamma^*$ is a lattice with
  $\beta_0\in\Gamma^*$ such that $\Gamma^*_{\beta_0}\not=\varnothing$
  then there exists a $\HP^n$--family of non congruent
  Hamiltonian stationary tori with Lagrangian angle
  $\beta=2\pi\langle\beta_0, \rangle$ given by
\begin{equation}
\label{eq:all hsl tori}
f=
 \ejbh \sum_{\delta\in\Gamma^*_{\beta_0, +}} \left(1-k\frac{2\delta}{\beta_0}\right) e_\delta c _\delta
\end{equation}
where  $c_\delta\in\H$, $e_\delta(z)=\edelta$,  $n+1=\dim_\H H^0(\trivial{})$ and
\[
\Gamma^*_{\beta_0,+} \assign \{\delta\in\Gamma^*_{\beta_0} \mid \Im(\delta\beta_0^{-1})>0\}\,.
\] 
Conversely, every Hamiltonian stationary torus with
Lagrangian angle $\beta_0$ arises, up to translation,
this way.
\end{theorem}
\begin{proof}
  We have seen that a given Hamiltonian stationary torus $f:
  T^2\to\R^4$ with Lagrangian angle $\beta=2\pi\langle
  \beta_0,\cdot\rangle$ has $\Gamma^*_{\beta_0}\not=\varnothing$ and
  $f\in H^0(\trivial{})$, that is $f$ is given, up to translation, by
  a complex linear combination of the holomorphic sections
  $\alpha_\delta$, $\delta\in\Gamma^*_{\beta_0}$, with trivial
  multiplier. Corollary \ref{cor:holomorphic sections} thus gives 
\[
f =
  \ejbh \sum_{\delta\in\Gamma^*_{\beta_0}}
  \left(1-k\frac{2\delta}{\beta_0}\right) e_\delta \tilde u_\delta
\]
with $\tilde u_\delta\in\C$. Since $\Gamma^*_{\beta_0}$ is symmetric,
in particular $H^0(\trivial{})$ is quaternionic, we can use only half
of the lattice points in $\Gamma^*_{\beta_0}$ by replacing the complex
Fourier coefficients by quaternionic ones: using
$\frac{|\beta_0|}{2|\delta|}=1$ for $\delta\in\Gamma^*_{\beta_0}$ and
$-k e_\delta k = e_{-\delta}$ we have
\[
\left(1-k \frac{2\delta}{\beta_0}\right) e_\delta \tilde u_\delta +
\left(1+ k \frac{2\delta}{\beta_0}\right) e_{-\delta} \tilde
u_{-\delta} = \left(1-k \frac{2\delta}{\beta_0}\right)
e_\delta\left(\tilde u_\delta + \frac{\beta_0}{2\delta} k \tilde
  u_{-\delta}\right)
\] 
which immediately gives (\ref{eq:all hsl tori}) with $c_\delta=
(\tilde u_\delta + \frac{\beta_0}{2\delta} k \tilde
u_{-\delta})\in\H$. 

  Conversely, given a lattice $\Gamma$ and a Lagrangian angle
  $\beta=2\pi\langle\beta_0, \cdot \rangle, \beta_0\in\Gamma^*$ with
  $\Gamma^*_{\beta_0}\not=\varnothing$, we can define a complex
  structure on the trivial bundle $T^2\times\H$ over the torus
  $T^2=\C/\Gamma$ by left multiplication by $N = e^{j\beta}i$, and a
  holomorphic structure by $D = d''$. In particular, every holomorphic
  section of $D$ gives a conformal torus $f: T^2\to\H$ with left
  normal $N=e^{j\beta}i$, and thus $f$ is Hamiltonian stationary.  Since every holomorphic section with trivial multiplier is
given by (\ref{eq:all hsl tori}) and  
\[
n+1 =\dim_\H H^0(\trivial{})= |\Gamma^*_{\beta_0,+}| + 1
\]
we obtain this way a $\HP^n$--family of non congruent
Hamiltonian stationary tori.

\end{proof}

\begin{rem} We have recovered the result of \cite{pascal_quaternionic}
  without integration of $df = \ejbh dz g$. Moreover, the previous
  theorem shows that the family of Hamiltonian stationary tori with
  the same lattice $\Gamma$ and Lagrangian angle $\beta$ is obtained
  by projections of a quaternionic holomorphic curve $F:
  \C/\Gamma\to\HP^k$ to $\HP^1$ whose complex structure, cf. for
  example \cite{habil}, is given by $N=e^{j\beta} i$.
\end{rem}

\begin{theorem}
\label{thm:real structure}
The real structure $\rho: \Spec \to \Spec$ induces involutions 
\[
\rho: \C^2\to\C^2, (A,B) \mapsto (A, -B)
\quad \text{
and
} \quad 
\rho: \Gamma^* \to \Gamma^*, \delta \mapsto -\delta
\]
so that $ \rho(\Gamma^*_{A,B}) = \Gamma^*_{A,-B}$ and  
$\L$ is compatible with $\rho$ that is
\[
\rho^* \L_{A,\delta-B} \assign \L_{A,-(\delta-B)} = \L_{A,\delta-B} \, j\,.
\]

\end{theorem}
\begin{proof}
  Note that $\rho(h^{A,B}) = h^{A,-B}$, and from (\ref{eq: Gamma 0B})
  and (\ref{eq:Gamma_A,B, else}) we see that $\delta\in\Gamma^*_{A,B}$
  implies $-\delta\in\Gamma^*_{A,-B}$.  In Theorem
  \ref{thm:parametrized multipliers} the complex line
  $\L_{A,\delta-B}$ is defined as the span of the monochromatic
  holomorphic section
\[
\alpha^{A,B}_\delta = \ejbh(1-k \lambda_{A,\delta-B})e_{\delta-B} e^{2\pi \langle A, \cdot \rangle}\,.
\]
 By (\ref{eq:lambda delta}), (\ref{ABcondition}) we see
$
\lambda_{A,\delta-B} \bar\lambda_{A, -(\delta-B)} %= \frac{4}{|\beta_0|^2}(-(\bar\delta-\bar B) + i\bar A)((\delta-B) -iA) 
= -1
$
so that
\[
\alpha_{-\delta}^{A,-B}  j = \alpha_\delta^{A,B} w
\]
with $w = i\bar\lambda_{A,-(\delta-B)} \in\C$, and $\L$ is
compatible with $\rho$.
\end{proof}
A real multiplier $h\in\Spec$, that is
$h_\gamma = \bar h_\gamma$ for all $\gamma\in\Gamma^*$, has a
quaternionic space of holomorphic sections since $\alpha\in
H^0_h(\ttrivial)$ implies $\alpha j\in H^0_h(\ttrivial)$. In particular, 
\[
\dim_\C H^0_h(\ttrivial) = 2\dim_\H
H^0_h(\ttrivial)\ge 2\,,
\]
and in fact, real multipliers occur at the double points $A$:
 
\begin{cor}
A multiplier $h\in\Spec$ is real if and only if $h=h^{A,B}$ for a double point 
\[
 A = -\frac{i}2
  \zeta \sqrt{1 - \frac{|\beta_0|^2}{|\zeta|^2}}, \quad \zeta\in\Gamma^*, \quad |\zeta| > |\beta_0|\,,
\]
 or $h=h^{0,B}$ with  $B\in\frac 12\Gamma^*$.
\end{cor}
\begin{proof}

  If $h = h^{0,B}$ then $B =\frac{\beta_0}2(1+ e^{it_0})\mod\Gamma^*$
  for some $t_0\in[0,2\pi)$. But then
\[
h^{0,B} = e^{-\pi i
    \langle\beta_0(1+e^{it_0}), \cdot\rangle}
\]
is real if and only if $\beta_0(1+e^{it_0}) \in\Gamma^*$, that is
$B\in\frac 12\Gamma^*$ in which case $h^{0,B}_\gamma = e^{2\pi
  i\langle B, \gamma\rangle} =\pm 1$ for all $\gamma\in\Gamma^*$.

If $A\not=0$ and $h=h^{A,B}$ is real then $\dim_\C H^0_h(\ttrivial)\ge
2$ shows that $A$ is a double point.  Conversely, if 
\[
A = -\frac
i2\zeta \sqrt{1-\frac{|\beta_0|^2}{|\zeta|^2}}, \quad \zeta\in\Gamma^*,
\quad |\zeta|>|\beta_0|
\]
is a double point then we obtain from (\ref{eq:Gamma_A,B, else}) that  $B= \frac 12(\beta_0-\zeta)\hspace{-6pt}\mod\Gamma^*$ which shows that $B\in\frac 12\Gamma^*$ and 
$
h^{A,B}_\gamma = \pm e^{2\pi\langle A, \gamma\rangle} \in\R$ is real.
\end{proof}

\begin{cor}
  Over generic points of the spectral curve $\Sigma$ the kernel
  bundle $\L$ is given by
\[
\mathcal L_h = \L_{A,\delta-B}\,,
\]
and $\L$ is compatible with the real structure $\rho$, that is
$\rho^*\L = \L j$.
\end{cor}

Summarizing the previous results we see that for a generic multiplier
$h=h^{A,B}$ there is a 1--dimensional space of holomorphic
sections. If we denote 
\[
\Spec_0 =\{ h \in\Spec \mid \dim
H^0_h(\ttrivial)=1\}
\] we have a well--defined map
\[
\lambda: \Spec_0 \to\C_*, \quad h \mapsto \lambda(h) = \lambda_{A,\delta-B}, 
\]
where the multiplier $h=h^{A,B}$ has $\Gamma^*_{A,B} =\{ \delta\}$ and
$\lambda_{A, \delta-B}$ is defined by (\ref{eq:lambda delta})
\begin{equation}
\label{eq:lambda_again}
\lambda_{A,\delta-B} = -\frac{2i(|A|\mp r)A}{\beta_0|A|}, \quad 
r= \sqrt{|A|^2 + \frac{|\beta_0|^2}4} \,.
\end{equation}
Here
we used that the expression $\delta-B$ is uniquely defined by $h$
since $\Gamma^*_{A,B}= \Gamma^*_{A,B+\zeta} -\zeta$ for all $\zeta\in
\Gamma^*$.  The map $\lambda$ gives rise to the normalization of
$\Spec$:
\begin{prop}
\label{prop:normalization}
There is a surjective map $\eta : \mathbbm{C}_{\ast} \rightarrow
\tmop{Spec}$ with $\eta(\lambda(h)) = h$ for all $h\in\Spec_0$.  The
pullback bundle $\eta^{\ast} \mathcal{L}$ extends smoothly across
$\lambda (\Spec \setminus \Spec_0) \subset \mathbbm{C}_{\ast}$.
\end{prop}
\begin{proof} 
  For $\lambda\in\C_*$ we define $\eta(\lambda)\assign
  h^{A_\lambda,B_\lambda}$ where $A_\lambda = \frac{i \beta_0}{4}
  (\lambda - \bar{\lambda}^{- 1})$ and $B_\lambda = \frac{\beta_0}{4}
  (2-\lambda\invers - \bar{\lambda}).$ For $h\in\Spec_0$, $h=h^{A,B}$,
  with $\lambda(h)\not\in S^1$ we see from (\ref{eq:lambda_again}) that 
\[
A = \frac{|\lambda(h)|^2-1}{4\overline{\lambda(h)}}\,.  
\]
Clearly, this extends to $h\in\Spec_0$ with $\lambda(h)\in S^1$, that
is $h=h^{A,B}$ with $A=0$, and thus $A=A_\lambda$ for all $h\in
\Spec_0$, $\lambda=\lambda(h)$. Using (\ref{eq:lambda delta})   we see
\[
B = \delta - \frac{\beta_0}2\lambda(h) + i A = \delta - \frac{\beta_0}4(\lambda(h) + \overline{\lambda(h)}\invers)
\]
with $\delta\in\Gamma^*_{A,B} \subset \Gamma^* + \frac{\beta_0}2$ so
that $B=B_\lambda \mod \Gamma^*$, and $\eta(\lambda(h^{A,B})) = h^{A,
  B}$.  Finally, Corollary \ref{cor:holomorphic sections} shows that
the line bundle $\eta^{\ast} \mathcal{L}$ is well defined over $\C_*$.
\end{proof}

We will see later that $\mu=\frac{1}{\lambda^2}$ is in fact
(\ref{eq:spec_parameter}) the spectral parameter of the harmonic map
given by the left normal $N$ of $f$.

%%% Local Variables: 
%%% mode: latex
%%% TeX-master: "doc"
%%% End: 

%% file: darboux.tex
\section{Darboux transforms}
\label{sec:Darboux}

In the previous section, we discussed all possible multipliers of
holomorphic sections of the quaternionic holomorphic line bundle which
is associated to a Hamiltonian stationary torus. Each holomorphic
section with multiplier gives rise \cite{conformal_tori} to a
conformal map, a Darboux transform of $f$. We quickly recall this
construction, and show that the Darboux transforms of a Hamiltonian
stationary torus which arise from monochromatic holomorphic sections
are again Hamiltonian stationary.

Recall that a sphere congruence $S\colon M \to S^4$ assigns to each
point $p\in M$ a 2--sphere $S(p)$. If the sphere congruence $S$ passes
through a conformal immersion $f\colon M \to\R^4$, that is $f(p)\in
S(p)$, and the Gauss maps of $f$ and $S(p)$ coincide at each point
$p\in M$, then $S$ is said to \emph{envelope} $f$. A pair of conformal
immersions $f, f^\sharp\colon M \to\R^4$ such that there exists a
sphere congruence $S$ enveloping both $f$ and $f^\sharp$ is called a
\emph{classical Darboux pair}. In this case, both $f$ and $f^\sharp$
are isothermic \cite{darboux}. More generally, a branched conformal
immersion $\hat f: M \to S^4$ is called a Darboux transform of a
conformal immersion $f: M \to S^4$ if there exists a sphere congruence
$S$ on $M$ enveloping $f$ and \emph{left--enveloping} $\hat f$ over
immersed points. The enveloping conditions can be written in terms of
complex structures on the trivial $\H^2$--bundle $V$ over $M$: $\hat
f$ is a Darboux transform of $f$ if and only if there exists a complex
structure $S\in\Gamma(\End(V))$ such that $L$ and $\hat L$ are
$S$--stable and $*\delta = S \delta = \delta S$ and $*\hat \delta =
S\hat\delta$ where $\delta$ and $\hat \delta$ denote the derivatives
of the line bundles $L$ and $\hat L$ given by $f$ and $\hat f$
respectively.  Note that the global existence of $S$ already implies
that $\hat f$ has only isolated branch points \cite{Klassiker}.  In
the case, when $f\colon M \to\R^4$ maps to the Euclidean 4--space a
sphere congruence $S$ passing through $f$ is left--enveloping $f$ if
the left--normals of $f$ and $S(p)$ coincide at each point $p\in M$.

\begin{lemma}[see \cite{conformal_tori}] Every holomorphic section $\varphi\in
  H^0(\widetilde{V/L})$  of the canonical
  holomorphic bundle of a conformal immersion $f: M \to S^4$ has a
  unique lift $\hat\varphi\in\Gamma(\widetilde{V})$ of $\varphi$ such
  that
\begin{equation}
\label{eq:prolongation}
\pi d\hat\varphi=0\,,
\end{equation}
where $\pi: V \to V/L$ is the canonical projection.  This unique lift
$\hat\varphi$ is called the \emph{prolongation} of $\varphi$.
\end{lemma}
Note that the prolongation $\hat\varphi$ has the same multiplier as
$\varphi$ so that, if $\varphi$ has no zeros, $\hat f =
\hat\varphi\H\colon M\to S^4$ defines a map from the Riemann surface
$M$ into the 4--sphere which turns out to be a 
 Darboux transform of $f$. In the case when $\varphi$ has
zeros, one obtains a conformal map $\hat f$ away from the zeros of
$\varphi$, which is again a Darboux transform on its domain.  Such a
map $\hat f$ is called a \emph{singular} Darboux transform.
Conversely, every (singular) Darboux transform $\hat f\colon M\to S^4$
gives a holomorphic section with multiplier of $V/L$. 
\begin{definition} Let $f: T^2\to\R^4$ be a Hamiltonian stationary
  torus with associated holomorphic line bundle $V/L$.

  A branched conformal immersion $\hat f: T^2 \to S^4$ is called a
  \emph{monochromatic Darboux transform of $f$} (respectively
  \emph{polychromatic}) if it is given by the prolongation of a
  monochromatic (respectively polychromatic) holomorphic section with
  multiplier of $\widetilde{V/L}$.
\end{definition}

For a conformal immersion  $f: M \to \R^4$ we identify as before $V/L =
\trivial{}$ via the nowhere vanishing section $\pi e\in\Gamma(V/L)$
where $e=\begin{pmatrix} 1 \\0
\end{pmatrix}
$.  To compute the prolongation in terms of the trivialization, note
that every nowhere vanishing holomorphic section $\alpha\in
H^0_h(\widetilde{\trivial{}})$ with multiplier $h$ defines a unique
quaternionic flat connection $\hat d$ on $\ttrivial$ by requiring
$\hat d \alpha=0$. Since $\alpha $ is holomorphic, that is $*d\alpha =
Nd\alpha$, we can write
\[
\hat d = d + df \hat T\,,
\]
where $\hat T: M \to \H$ is defined on $M$ rather than on the universal 
cover $\tilde M$ of $M$ since $\alpha$ and $d\alpha$ have the same
multiplier.

\begin{lemma}
\label{lem:prolongation}
Let $f: M \to\R^4$ be a conformal immersion.  The prolongation $\hat
\varphi$ of a nowhere vanishing, non--constant holomorphic section
$\varphi = \pi e\alpha\in H^0(\widetilde{V/L})$ is given by
\[
\hat \varphi = \psi \nu + e \alpha
\]
where $\psi =\begin{pmatrix} f\\1
\end{pmatrix}
\in\Gamma(L)$, $\nu = 
\hat T \alpha$, and  $\hat T$ is defined by $df \hat T \alpha =-d\alpha$.\\

In particular, if the Darboux transform $\hat f: M \to \R^4$
associated to $\alpha$ maps into $\R^4$ rather than $S^4$, that is, if
$\hat T$ is nowhere vanishing, then $ \hat f = f + T: M \to\R^4$ with
$T=\hat T\invers$.
\end{lemma}
\begin{proof}
  One easily verifies that $\hat \varphi$ satisfies
  (\ref{eq:prolongation}). But then, if $\hat T\not=0$, 
\[
\hat f = \hat\varphi \H = \begin{pmatrix} f \nu + \alpha\\ \nu
\end{pmatrix}
 \H = \begin{pmatrix} f + T \\1
 \end{pmatrix}
 \H
\]
and $\hat f = f +T$ is the Darboux transform of $f$ when choosing the
point at infinity $\infty =e\H$.
\end{proof} 
\begin{rem} The previous lemma extends in the obvious way to
  holomorphic sections with zeros and singular Darboux transforms.
\end{rem}
From (\ref{eq:holomorphicity}) we know that $\Span_\H\{1,f\}\subset
H^0(\trivial{})$. In particular, if $\alpha\in H^0(\trivial{})$ is
constant the previous lemma shows that $\hat\varphi = e\alpha$ is the
prolongation of $\varphi =\pi e\alpha$.  In other words, the Darboux
transform given by $\alpha$ is the constant map $\hat f =
\infty$. More generally,
\begin{cor}
\label{cor:trivial multiplier}
If $f: T^2\to\R^4$ is a conformal immersion and $\alpha\in
\Span_\H\{1,f\}\subset H^0(\trivial{})$, then the Darboux transform of
$f$ given by $\alpha$ is a constant map $\hat f = c$ with
$c\in\H\cup\{\infty\}$.

  In particular, if the complex dimension of the space of global
  holomorphic sections is minimal, that is $\dim_\C
  H^0(\trivial{})=4$, then every Darboux transform given by the
  trivial multiplier is constant.
\end{cor}
 
From (\ref{eq:monochromatic holomorphic section}) we see that
monochromatic holomorphic sections of a Hamiltonian stationary torus
are always nowhere vanishing, in particular, we obtain regular Darboux
transforms defined on the torus.  Lemma \ref{lem:prolongation} allows
to compute all Darboux transforms in 4--space of a Hamiltonian
stationary torus:

\begin{theorem}
\label{thm:Darboux transforms}
Every non--constant Darboux transform $\hat f: T^2\to\R^4$ of a
Hamiltonian stationary torus $f$ with $df = \ejbh dz g$ is given by a
nowhere vanishing holomorphic section $\alpha\in H^0_h(\ttrivial)$
with multiplier $h=h^{A,B}$, and
\begin{equation}
\label{eq: monochromatic tau, else}
     \hat f=       f - \ejbh \frac{1}{\pi(4|A|^4 +\langle\beta_0, A\rangle^2)}(2|A|^2 - j \langle\beta_0, A\rangle)A g\,,
       \quad  \text{ if } A\not =0\,,
\end{equation}       
or 
\begin{equation}
\label{eq: multiple darboux}
\hat f = f + \ejbh \left(\sum_{s,t\in I_B} (1+k e^{is}) u_s \bar u_t e_{\frac{\beta_0}2(e^{it} -e^{is})} (1 + ke^{it})  \sin t\right) \frac{1}{r \pi \bar\beta_0} g
\end{equation}
where   the finite set
\[
I_B=\{ t\in [0, 2\pi) \mid B-\frac{\beta_0}2e^{it}\in\Gamma^*_{0,B}\}\not=\{0,\pi\}
\]
parametrizes the admissible frequencies and $u_t\in \C$ are chosen so
that the map
\[
r =|\sum_{t\in I_B} u_t \sin t \,
  e_{B-\frac{\beta_0}2e^{it}}|^2 +
|\sum_{t\in I_B} u_t e^{it}\sin t\, 
  e_{B-\frac{\beta_0}2e^{it}}|^2 
\]
is nowhere vanishing.  Here we again use $e_\delta(z) = \edelta$ for
$\delta\in\Gamma^*+\frac{\beta_0}2$. Equation \eqref{eq: multiple
  darboux} simplifies in the case of a monochromatic holomorphic
section to
\begin{equation}
\label{eq: monochromatic tau, A=0}
\hat f=  f + \ejbh k\frac{e^{it}}{\pi\bar\beta_0\sin(t)} g\,,\quad  \text{ if }  A=0, I_B = \{t\}, t\not\in\{0,\pi\}\,.
\end{equation}
\end{theorem}
\begin{rem}
  The compatibility of the real structure $\rho$ with
  $\L_{A,\delta-B}$ implies, see Theorem \ref{thm:real structure},
  that the induced Darboux transforms of $h^{A,B}$ and $\rho(h^{A,B})
  = h^{A,-B}$ coincide. Therefore, in the case when $A\not=0$ the
  Darboux transform does not depend on whether $\Gamma^*_{A,B} =
  \{\delta_+\}$ or $\Gamma^*_{A,B} = \{\delta_-\}$.  In particular, if
  $h=h^{A,B}$ is a multiplier at a double point $A$, then $h$ is real,
  and we obtain only one closed Darboux transform, given by (\ref{eq:
    monochromatic tau, else}), associated with the multiplier $h$.
\end{rem}

\begin{proof}[Proof of Theorem \ref{thm:Darboux transforms}]
  It is shown in \cite{conformal_tori} that every Darboux transform is
  given by a holomorphic section $\alpha\in H^0_h(\ttrivial)$ and we
  have seen that $h=h^{A,B}$.  To compute globally defined Darboux
  transforms $\hat f: T^2\to\R^4$ with values in Euclidean 4--space, we
  have to find holomorphic sections $\alpha$ with multiplier so that
  the flat quaternionic connection $\hat d = d + df \hat T$ with $\hat
  d \alpha=0$ satisfies $\hat T\not=0$.  As before it is
  computationally easier to use the gauged holomorphicity condition
  for $\tilde \alpha = \emjbh \alpha$: the unique flat quaternionic
  connection with $\tilde d\tilde \alpha=0$ is
\[
\tilde d = d + \frac{jd\beta}2 + dz \hat \tau
\]
with $\hat \tau=g \hat T\ejbh$; in particular, it is enough to find
$\tilde \alpha$ with nowhere vanishing $\hat \tau$. By Corollary
\ref{cor:holomorphic sections} all holomorphic sections $\alpha$ with
multiplier $h=h^{A,B}$ are given by
\[
\tilde\alpha = \tilde \sigma e^{2\pi(\langle A, \cdot\rangle - i \langle B, \cdot \rangle)}\,
\]
where  $\tilde \sigma
=\sum_{\delta\in\Gamma^*_{A,B}}(1-k\lambda_\delta)e_{\delta}\tilde
u_\delta$ so that
\begin{eqnarray*}
 \lefteqn{(\tilde d \tilde \alpha) e ^{-2\pi(\langle A, \cdot\rangle - i \langle B, \cdot \rangle)}}\\ 
 &=&\pi\sum_{\delta\in\Gamma^*_{A,B}} \Big[(1 -k  \lambda_\delta)\Big((i\bar\delta - i \bar B + \bar A)dz + (i\delta- i B  + A)d\bar z\Big) \\
&& \hspace{2cm}  + \frac{j}2(\bar\beta_0 dz + \beta_0 d\bar z)(1 -k \lambda_\delta)\Big] \tilde u_\delta e_\delta 
+ dz \hat\tau \tilde \sigma\\
&=& \pi d\bar z \sum_{\delta\in\Gamma^*_{A,B}}\left( (i\delta -i B +  A) - k \lambda_\delta (i\bar\delta - i\bar B + \bar A) + \frac{j}2 \bar\beta_0(1 -k\lambda_\delta)\right) \tilde u_\delta e_\delta \\
&& +  d z \left( \pi\sum_{\delta\in\Gamma^*_{A,B}}\left( (i\bar\delta -i \bar B + \bar A) - k\lambda_\delta (i\delta - i B + A) + \frac{j}2 \beta_0(1 -k\lambda_\delta)\right) \tilde u_\delta e_\delta + \hat\tau \tilde \sigma\right)\,.
\end{eqnarray*}
By holomorphicity of $\alpha$ the $(0,1)$--part  with respect to left multiplication by $i$
\[
 \pi d\bar z\sum_{\delta\in\Gamma^*_{A,B}}\left( (i\delta -i B +  A) + j i \lambda_\delta (i\bar\delta - i\bar B + \bar A) + \frac{j}2 \bar\beta_0(1 + ji\lambda_\delta)\right) \tilde u_\delta e_\delta =0
\]
vanishes (\ref{ABcondition}),  and thus $\hat \tau$ with $-d\tilde \alpha= \frac{jd\beta}2\tilde \alpha + dz \hat \tau\tilde \alpha$ satisfies
\begin{equation}
\label{eq:general_tau}
\hat\tau\tilde \sigma = - \pi \sum_{\delta\in\Gamma^*_{A,B}}\left( (i\bar\delta -i \bar B + \bar A) + (i\bar \delta - i \bar B - \bar A) k  \lambda_\delta + \frac{j}2 \beta_0(1 - k\lambda_\delta)\right) \tilde u_\delta e_\delta\,.
\end{equation}
If $\hat f$ is given by a monochromatic holomorphic section then
$\tilde \sigma = (1-k\lambda_\delta) e_\delta\tilde u_\delta$ is given
by a single frequency $\delta\in\Gamma^*_{A,B}$ and
\[
\hat \tau = -\pi\left\{\Big(i(\bar\delta -\bar B)\frac{1-|\lambda_\delta|^2}{1+|\lambda_\delta|^2} + \bar A\Big) + j\left(\frac{\beta_0}2 - 2(\delta-B)\frac{\lambda_\delta}{1 + |\lambda_\delta|^2}\right)\right\}
\]
is constant. Moreover,   (\ref{eq:lambda delta})
yields
\[
1 -|\lambda_\delta|^2 = \begin{cases} 0 & A=0\\
                                     \frac{8|A|}{|\beta_0|^2}(-|A| \pm r) & A \not=0\,, \, \delta =  B \pm i r\frac{A}{|A|}
                                   \end{cases}                                 
\]
and
\[
1 +|\lambda_\delta|^2 = \begin{cases} 2 & A=0\\
                                     \frac{8}{|\beta_0|^2}r(r \mp |A|) & A \not=0\,, \, \delta =  B \pm i r\frac{A}{|A|} 
                                   \end{cases} \,,
\] 
so that
\begin{equation*}
\hat \tau = \begin{cases} -k \pi\beta_0e^{it} 
\sin t  \quad & \text{ if }  A=0, \, \, \delta=B -\frac{\beta_0}2e^{it}\\[.2cm]
                      -\pi\left(2\bar A + j\frac{\langle \beta_0, A\rangle}{\bar A}\right) \quad & \text{ if } A\not =0
                    \end{cases}\,.
\end{equation*}  
If $A\not=0$ or if $A=0$ and $t\not\in\{0,\pi\}$ then $\hat\tau \not=
0$ is a non zero constant, and we obtain the Darboux transforms $\hat
f = f + \ejbh \hat\tau\invers g:T^2\to\R^4$ given by (\ref{eq:
  monochromatic tau, A=0}) and (\ref{eq: monochromatic tau, else}).

Finally, the only Darboux transforms which arise in families are given
by multipliers $h=h^{0,B}$ with $\dim H^0_{0,B} >2$ or $\dim H^0_{0,B}
= 2$ and $h$ not real.  In this case, (\ref{eq:general_tau})
simplifies to
\[
\hat\tau
=\pi\bar\beta_0 \left(\sum_{t\in I_B} (1  - k  e^{it}) \sin t  \,   u_t e_{B-\frac{\beta_0}2e^{it}} \right)\tilde\sigma\invers
\]
with $\tilde\sigma = \sum_{t\in I_B} (1 + ke^{it}) u_t
e_{B-\frac{\beta_0}2e^{it}}$, $u_t\in\C$, and we obtain the general
formula for $h=h^{0,B}$.
\end{proof}
 
If $f: T^2\to\R^4$ be a Hamiltonian stationary torus with Lagrangian
angle $\beta$ and $df=\ejbh dz g$ then the left normal of a Darboux
transform $\hat f = f + \ejbh \tau g$ of $f$ can be expressed entirely
in terms of $\tau$ and the Lagrangian angle $\beta$ of $f$: Lemma
\ref{lem:prolongation} and Theorem \ref{thm:Darboux transforms} show
that $d\alpha = -df T\invers \alpha$ for some holomorphic section
$\alpha\in H^0_h(\ttrivial)$ with multiplier $h$ where $T=\ejbh \tau
g$. Putting $\nu=T\invers\alpha$ this equation implies that $0 = df
\wedge d\nu$ and thus $*d\nu = - Rd\nu$ where $R$ is the right normal
of $f$. Therefore,
\[
d\hat f = df + dT = df + d\alpha \nu\invers - \alpha\nu\invers d\nu \nu\invers = - T d\nu \alpha\invers T
\]
shows that the left normal of $\hat f$ is given by $\hat N = - T R
T\invers$ since $*d\hat f = -TRT\invers d\hat f$.  But the right
normal of $f$ is given (\ref{eq:right normal}) by $R =-g\invers i g$
so that we obtain
\begin{equation}
\label{eq:left normal darboux}
\hat N = -TRT\invers = -\ejbh \tau i \tau\invers \emjbh\,.
\end{equation}

If $\hat f$ is a monochromatic Darboux transform, the left normal
$\hat N$ of $\hat f$ turns out to be harmonic:
\begin{theorem}
\label{thm:Darboux is hsl}
Let $f: T^2\to \R^4$ be a Hamiltonian stationary torus. Then every
Darboux transform of $f$, which is given by a holomorphic section in
$H^0_h(\ttrivial)$ with generic multiplier $h\in\Sigma$, is again
Hamiltonian stationary.

More precisely, if $\beta$ is the Lagrangian angle of $f$ then every
monochromatic Darboux transform $\hat f: T^2 \to \R^4$ of $f$ has
Lagrangian angle $\hat\beta= \beta + \beta_h$ where $\beta_h\in\R$ is
constant.
\end{theorem}
\begin{proof} 
If $\hat f$ is a monochromatic Darboux transform of $f$ then Theorem
\ref{thm:Darboux transforms} shows that $\tau = (\tau_0 + j\tau_1) c$
with $\tau_0, \tau_1\in\R, c\in \C$. In this case, the equation
(\ref{eq:left normal darboux}) for the left normal $\hat N$ of $ \hat
f$ simplifies to
\[
\hat N = \frac{1}{\tau_0^2 + \tau_1^2}\ejbh (\tau_0 + j\tau_1) i (\tau_0 -j\tau_1) \emjbh = e^{j\beta} \frac{(\tau_0+j\tau_1)^2}{\tau_0^2 + \tau_1^2} i\,.
\]
Since $ \frac{(\tau_0+j\tau_1)^2}{\tau_0^2 + \tau_1^2} \in S^3 \cap \Cj$, $\Cj =\Span_\R\{1, j\}$, we can write
\[
 \frac{(\tau_0+j\tau_1)^2}{\tau_0^2 + \tau_1^2} = e^{j \beta_h}
\]
with $\beta_h\in\R$ constant. This shows both that the left normal
$\hat N$ of the Darboux transform $\hat f$ takes values in the unit
circle $ S^3 \cap \Cj i$ and that $\hat N$ is harmonic, in other
words, $\hat f$ is Hamiltonian stationary.
\end{proof}

\begin{rem} In fact, we have seen that for all multiplier $h=h^{A,B}$
  with $A\not=0$ the corresponding Darboux transforms are Hamiltonian
  stationary. Moreover, for $A=0$ we may obtain Darboux transforms
  which are not Hamiltonian stationary only if $\dim_\C H^0_{0,B}>2$
  or $\dim_\C H^0_{0,B}=2$ and $h=h^{0,B}$ is not real.
\end{rem}

In Section \ref{sec:examples} we show that there exist Darboux
transforms which are not Lagrangian tori, and thus not Hamiltonian
stationary, in $\Cj^2$ using the following characterization of
Lagrangian Darboux transforms:
\begin{cor}
\label{cor:lagrangian}
Let $f: T^2\to\R^4$ be a Hamiltonian stationary torus.  A Darboux
transform $\hat f: T^2 \to \Cj^2$ of $f$ is Lagrangian in $\Cj^2$ if
and only $\hat f = f + \ejbh (\tau_0 + j\tau_1)g$ with $\Im(\bar\tau_0
\tau_1)=0$ where $\tau_0,\tau_1: M \to \C$.
\end{cor}
\begin{proof}
  From (\ref{eq:left normal darboux}) we see that the left normal
  $\hat N = n_0 i + n_1$ of $\hat f$ is given by
\[
n_0= \frac {|\tau_1|^2 - |\tau_0|^2 - 2 j \Re(\bar\tau_0\tau_1)}{|\tau_0|^2 + |\tau_1|^2} \, e^{j\beta}\in \Cj, 
\quad \text{ and } \quad  n_1 = 2j\frac{\Im(\bar\tau_0\tau_1)}{|\tau_0|^2 + |\tau_1|^2}  \in j \R\,.
\]
But then (\ref{eq:df}) shows that $\hat f$ is Lagrangian in $\Cj^2$ if
and only if $\hat N \in S^3\cap \Cj i$ which is equivalent to $n_1=0$.
\end{proof}

Given a Hamiltonian stationary torus $f$ with Lagrangian angle
$\beta$, Theorem \ref{thm:Darboux is hsl} shows that every
monochromatic Darboux transform $\hat f$ of $f$ has, after
reparametrization, again Lagrangian angle $\beta$.  Therefore, the
characterization of Hamiltonian stationary tori by there Lagrangian
angle and lattice in Theorem \ref{thm:hslLagrangianangle} shows that
every monochromatic Darboux transform of $f$ is, after
reparametrization, in the same family of Hamiltonian stationary tori
as $f$. In particular:

\begin{cor}
\label{cor:all monochromatic DT}
Let $f: \C/\Gamma\to\R^4$ be a Hamiltonian stationary torus.  If the
space of global holomorphic sections has complex dimension $\dim_\C
H^0(\trivial{})=4$ then every monochromatic Darboux transform of $f$
is, up to rotation, translation and reparametrization, $f$.
\end{cor}

%%% Local Variables: 
%%% mode: latex
%%% TeX-master: "doc"
%%% End: 

%% file: example.tex
 
\section{Examples}
\label{sec:examples}
In this section we illustrate the previous results in the case of
homogeneous tori and Castro--Urbano tori. Moreover, we shall see that
there are examples of Hamiltonian stationary tori which allow families
of Darboux transforms which are not Lagrangian surfaces in
$\Cj^2$. Each family is obtained from one of the finitely many
multipliers $h=h^{0,B}\in\Spec$ with $\dim_\C H^0(\ttrivial)>1$.

\subsection{Homogeneous tori}

We consider homogeneous tori $f: \C/\Gamma \to \C^2$,
\[
f(x,y) = \frac{1}{r_1} e^{2\pi j r_1 x} + i \frac 1{r_2} e^{2\pi j r_2 y}\,,
\]
 where the lattice is given by $\Gamma = \frac{1}{r_1} \Z \oplus \frac{i}{r_2} \Z$, $r_1, r_2>0$.
Since the derivative of $f$ can be written as
\[
df =  2\pi e^{\pi j(r_1 x- r_2 y)}\, dz\, j e^{\pi j(r_1 x+r_2 y)}
\]
the conformal immersion $f$ is Hamiltonian stationary with Lagrangian
angle
\[
\beta(z)= 2\pi(r_1 x-r_2 y),  \quad \text{ that is }  \quad \beta_0 = r_1-r_2 i\in\Gamma^*\,,
\]
and $df = \ejbh dz g$ with $ g= 2\pi j e^{\pi j(r_1 x+r_2 y)}\,.  $
Let us first discuss monochromatic Darboux transforms of homogeneous tori:
these are given by
\[
\hat f =  e^{2\pi j r_1x}\lambda_1
 + i e^{2\pi j r_2 y} \lambda_2\,,
\]
where we obtain for $A=A_0 + iA_1\not=0$ by (\ref{eq: monochromatic
  tau, else})
\[
\lambda_1 = \frac 1{r_1}-w A_0, \quad \lambda_2 =\frac 1{r_2} + \bar w A_1
\quad \text{
with} \quad w = \frac{2\langle\beta_0, A\rangle + 4|A|^2j}{4|A|^4 +
  \langle\beta_0,A\rangle^2}\,,
\]
and for $A=0, B=
\frac{\beta_0}2(1+e^{it})$ mod lattice $\Gamma^*$, by (\ref{eq:
  monochromatic tau, A=0})
\[
\lambda_1 =\frac{r_2^2-r_1^2 +2 r_1 r_2 \cot t}{r_1(r_1^2+r_2^2)}, \quad \lambda_2 = \frac{r_1^2-r_2^2 -2 r_1 r_2 \cot t}{r_2(r_1^2+r_2^2)}\,.
\]
Since $\hat f$ is conformal, we see that in both cases
$\frac{|\lambda_1|}{|\lambda_2|} = \frac{r_2}{r_1}$ so that $\hat f$
is homogeneous:
\begin{lemma}
  All monochromatic Darboux transforms of a homogeneous torus $f$ are
  again (after rescaling and reparametrization) the homogeneous torus
  $f$.
\end{lemma}

To consider polychromatic Darboux transforms we have to discuss the
points lying in the intersection $\mathcal{D} \cap (\Gamma^* +
\frac{\beta_0}2)$ of the disc of radius $|\beta_0|$ and the translated
dual lattice. Since this depends on the ratio of the radii $r_1$ and
$r_2$ we will here only discuss the existence of polychromatic Darboux
transforms and show that some of these polychromatic Darboux
transforms are not Lagrangian in $\Cj^2$. We will return to the
general case and the study of the geometric properties of
polychromatic Darboux transforms in a future paper.

\subsubsection{The Clifford torus}

The Clifford torus is given by $r_1=r_2=1$.  The points on the
translated dual lattice $\Gamma^*+\frac{1-i}2$ lying in the closed
disk $\mathcal D$ around $\frac{\beta_0}2= \frac{1-i}2$ of radius
$|\beta_0|= \sqrt{2}$ are
\[
\mathcal{E} =\frac{\beta_0}2 + \{0, \pm  1, \pm i, \pm (1+i), \pm (1-i)
\}
\,.
\]

\begin{figure}[h]
\includegraphics[height=7.5cm]{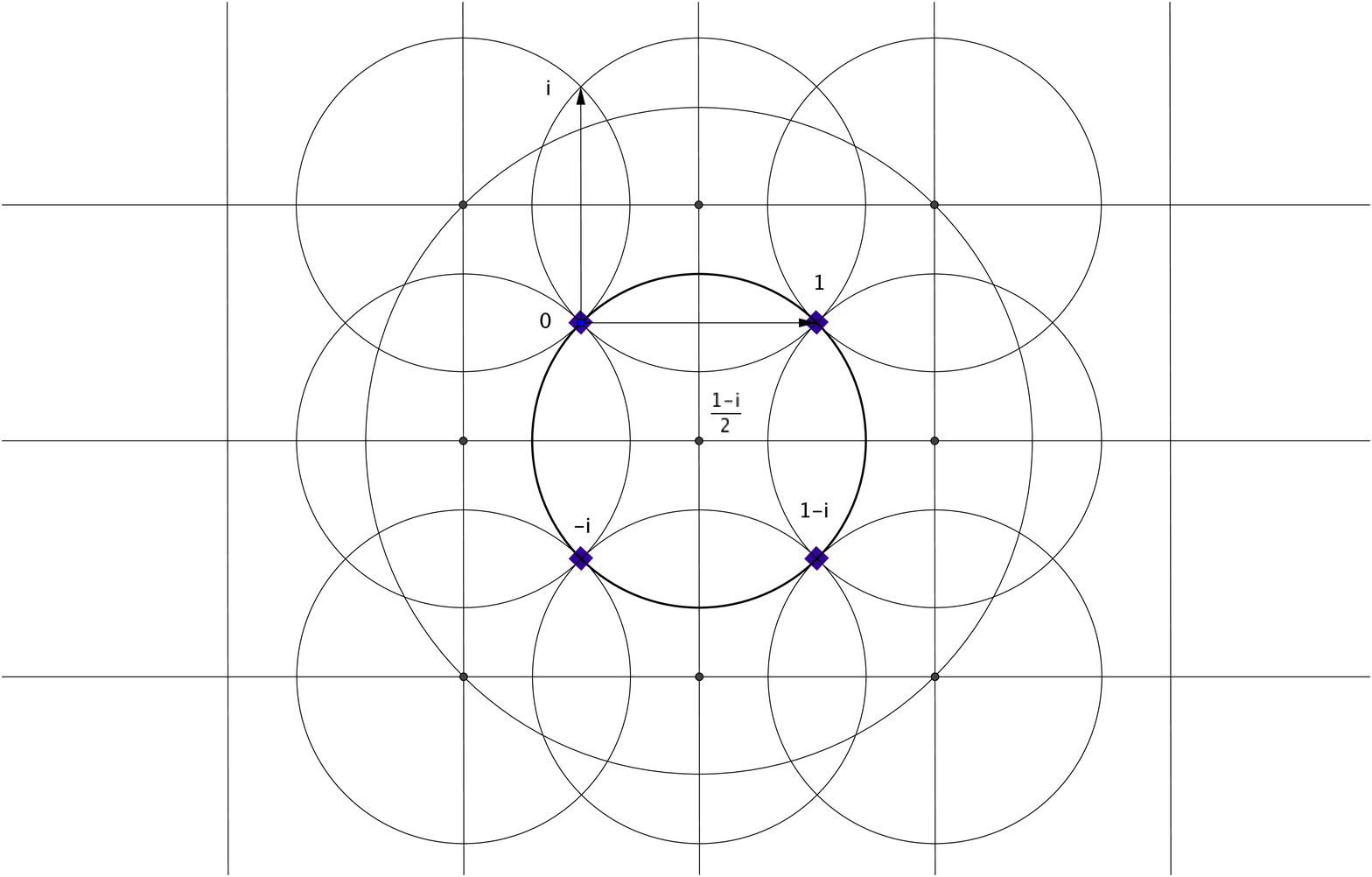}
\caption{}
\end{figure}

The circle $\mathcal{C}_{\frac{\sqrt{2}}2}(\epsilon)$ of radius
$\frac{\sqrt{2}}2$ around
$\epsilon\in\mathcal{E}\setminus\{\frac{\beta_0}2\}$ intersects the circle
$\mathcal{C}_{ \frac{\sqrt{2}}2}(\frac{\beta_0}2)$, and 
\[
 \mathcal{C}(\epsilon, \frac{\sqrt{2}}2) \cap
\mathcal{C}(\frac{\beta_0}2, \frac{\sqrt{2}}2) \subset\{ 0, 1, -i,
1-i\}\,.  
\]
Therefore, $\dim_\C H^0_{0, B} =1$ for all $B =
\frac{1-i}2(1 + e^{it}) $, $t\in(0,2\pi)\setminus\{\frac{\pi}2, \pi,
\frac{3\pi}2\}.$
 
The remaining points are $B = 0 \mod \Gamma^*$, that is $h = h^{0,0}$
is the trivial multiplier, and there are exactly four circles
intersecting at $B=0$, so that $ \dim H^0_{0, B} =4 $ for $B \in\{0,
1, -i, 1-i\}$.  In particular, the corresponding Darboux transforms
are by Corollary \ref{cor:trivial multiplier} constant maps, and there
are no polychromatic Darboux transforms of the Clifford torus:
\begin{lemma}
  Every non--constant Darboux transform $\hat f: T^2\to S^4$ of the
  Clifford torus is the (scaled and reparametrized) Clifford torus.
\end{lemma}
\begin{rem} Note that we only discuss closed Darboux transforms which
  are defined on the original torus $T^2 =\C/\Gamma$. If we allow
  other periods, then the examples of Bernstein \cite{holly} show that
  there are Darboux transforms (on an appropriate covering) of the
  Clifford torus which are isothermic but not constrained Willmore, in
  particular not Hamiltonian stationary.
\end{rem}

\subsubsection{The case $r_1=2$ and $r_2=1$}
The points in the translated dual lattice lying inside the closed disk
$\mathcal{D}$ around $\frac{\beta_0}{2} =\frac{2-i}{2}$ of radius
$|\beta_0|=\sqrt5$ are, see Figure 1, 
\[
\mathcal{E} =\frac{\beta_0}{2}+\{0, \pm 2, \pm i, \pm 2i , \pm (2+i), \pm (2-i)\}\,.
\]
The circle $\mathcal{C}_{\sqrt{5}/2}(\epsilon)$ of radius $\sqrt5/2$
around $\epsilon\in\mathcal{E}\backslash\{\frac{\beta_0}{2}\}$
intersects the circle $\mathcal{C}_{\sqrt{5}/2}(\frac{\beta_0}{2})$ in
at most two points, taken in $\frac{\beta_0}2 +\{1 \pm \frac i2, \frac 12 \pm i, -\frac 12 \pm i, -1\pm \frac i2 \} .  $ However each of
these points is congruent modulo $\Gamma^*$ to one of $0$,
$\frac{1+i}{2}$ or $\frac{3+i}2$.  The point $B=0$ lies at the
intersection of four circles, so $\dim_{\mathbb{C}} H^0_{0,0} = 4$,
and again the Darboux transforms given by the constant multiplier
$h=h^{0,0}$ are constant maps.

Let now $B=\frac{1+i}{2}$ (the case $B= \frac{3+i}2$ can be treated in
a completely analogous way). Then $\dim H^0_{0,B}=2$ since 
\[
\Gamma^*_{0,\frac{1+i}{2}} = \{
\frac{\beta_0}{2}, \frac{\beta_0}{2}+2i  \} \quad \text{ with } \quad I_B=\{\pi
-\arctan \frac{3}{4},  \frac{3\pi}{2}\}\,,
\]
and every polychromatic Darboux transforms, which is given by the
multiplier $h=h^{0,\frac{1+i}2}$, can be written (\ref{eq: multiple
  darboux}) as $\hat f = f + \ejbh \tau\frac{1}{r\pi\bar\beta_0} g$,
$\tau =\tau_0 + j\tau_1$ with
\begin{eqnarray*}
\tau_0 &=&
 (12i-16)u_1\bar u_2 e^{-4\pi i y} + 
20i\bar u_1 u_2 e^{4\pi i y}\\
\tau_1&=& (18 + 24i)|u_1|^2 + 50 |u_2|^2 
-(32i+24) u_1 \bar u_2 e^{-4\pi i y} 
+10\bar u_1 u_2 e^{4\pi i y}\,,
\end{eqnarray*}
where $u_1, u_2\in\C$ so that
\[
r=18|u_1|^2 + 50|u_2|^2 -12\, \Re\Big((1+2i)\bar u_1 u_2 e^{4\pi iy}\Big)
\]
is nowhere vanishing.  In particular, we see that in general
$\Im(\bar\tau_0\tau_1)\not\equiv 0$, so that Corollary
\ref{cor:lagrangian} shows that the family of polychromatic Darboux
transforms $\hat f$ given by the multiplier $h=h^{0,\frac{1+i}2}$ is
not Lagrangian, and thus not Hamiltonian stationary, in $\Cj^2$.

\subsubsection{The cases $r_1=m\in\N, m\ge 2$ and $r_2=1$}

Analogous to the case $r_1=2$, we see that $B=\frac{m-1}2(1+i)$ lies
at the intersection of two circles, and
\[
\Gamma^*_{0,\frac{m-1}2(1+i)}  = \{\frac{\beta_0}2, \frac{\beta_0}2+ mi\}
\]
so that we obtain a family of polychromatic Darboux transforms by
(\ref{eq: multiple darboux}) with $I_B = \left\{ 
  \pi-\arctan\left(\frac{m^2-1}{2m}\right), \frac{3\pi}{2}\right\}$.

\begin{figure}[h]
\includegraphics[height=4.5cm]{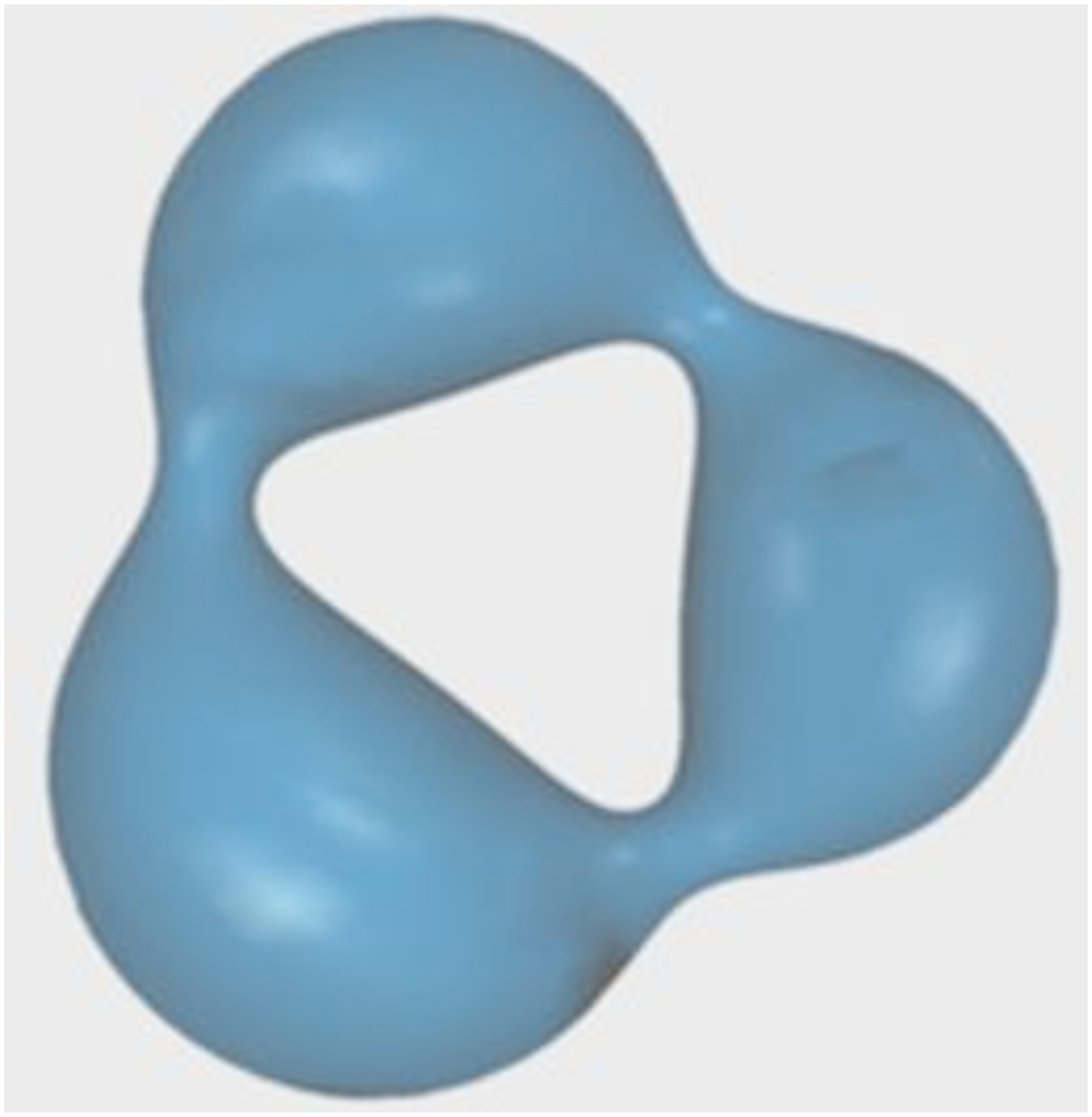} 
\includegraphics[height=4.5cm]{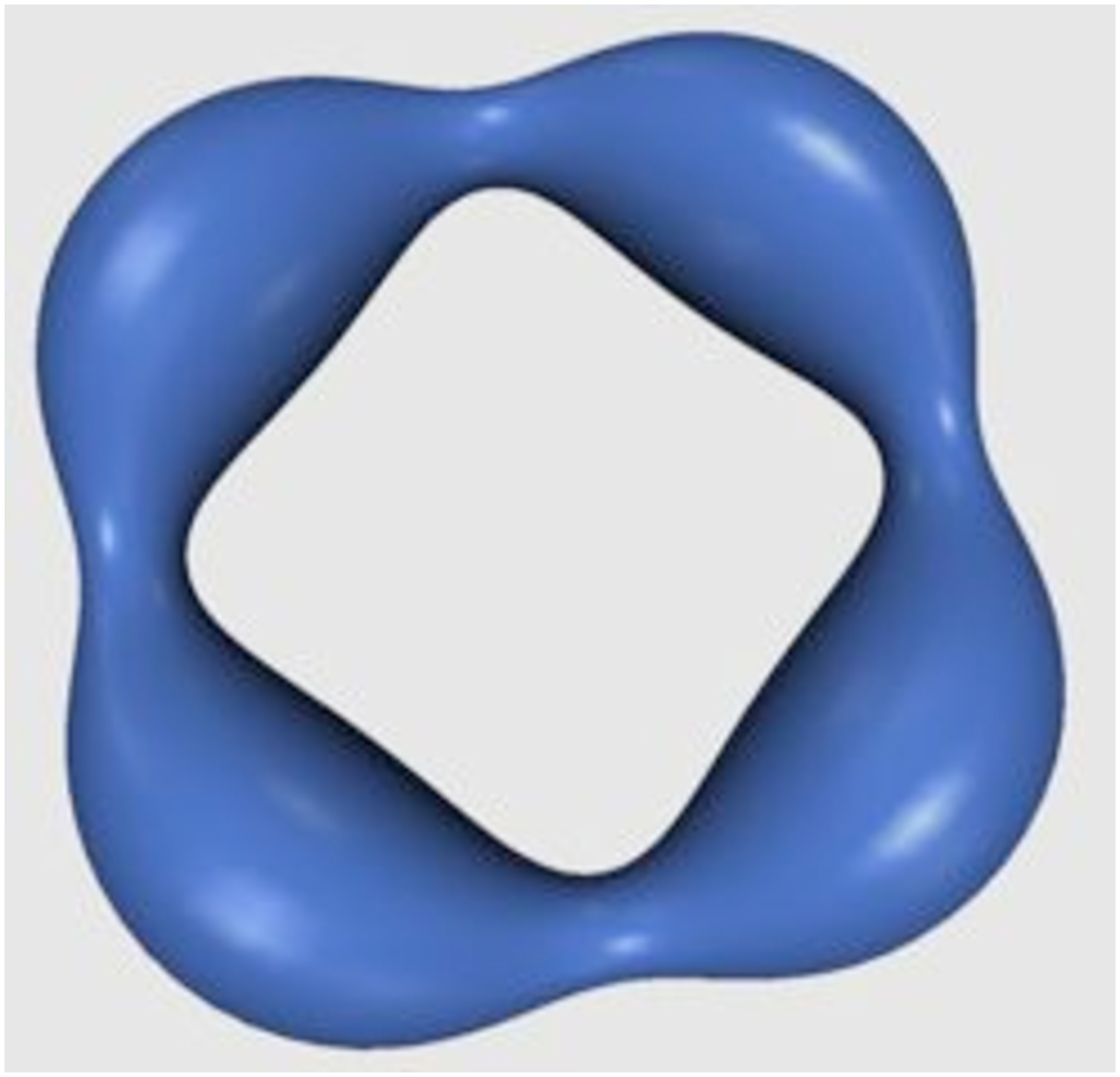}
\includegraphics[height=4.5cm]{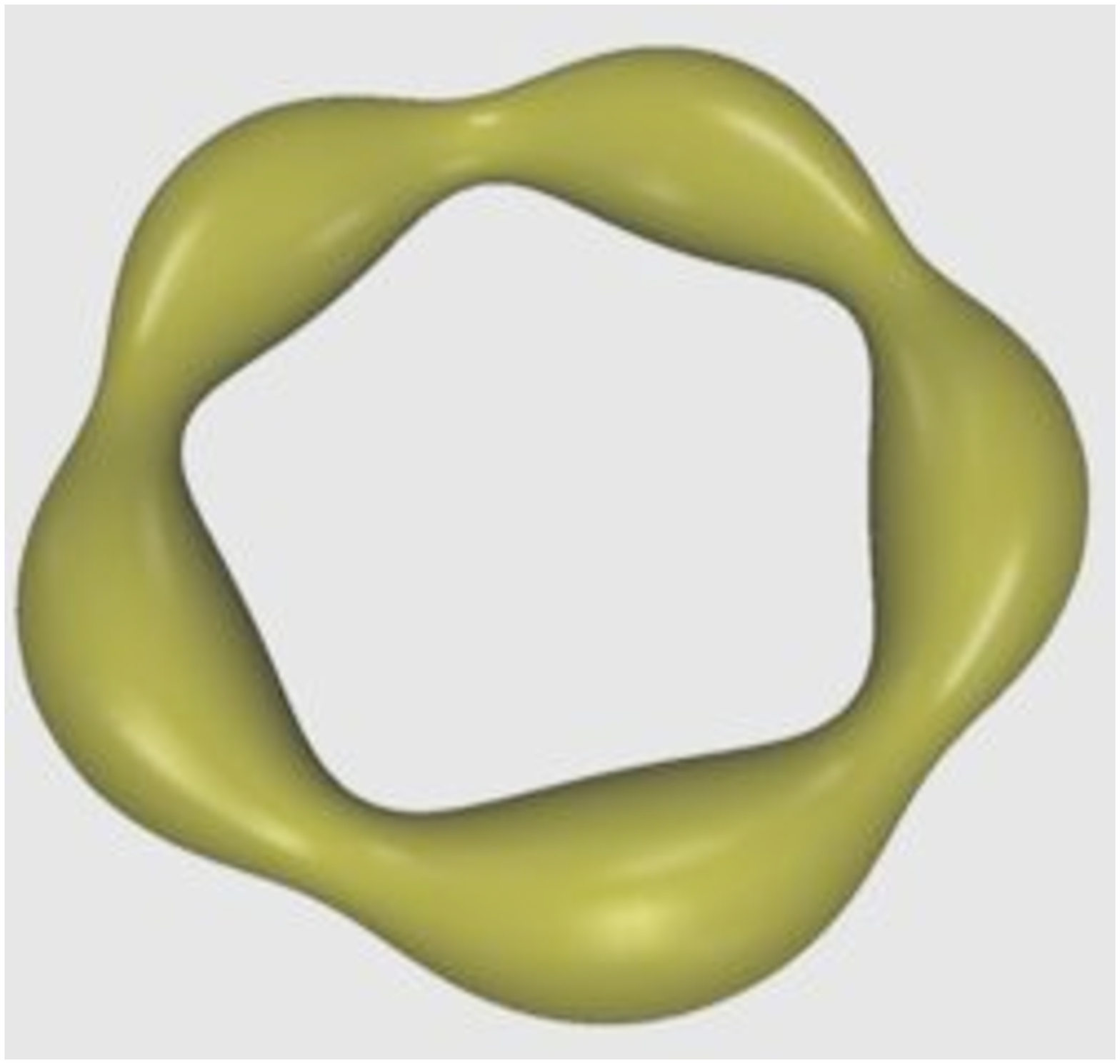}
\caption{Polychromatic Darboux transforms of homogeneous tori with
  $r_1=3,4,5$ and $r_2=1$}
\end{figure}

A similar computation as for the case $r_1=2$ shows that these
surfaces are again not Lagrangian surfaces in $\Cj^2$.  Note that
there may exist further $B$ with $|\Gamma^*_{0,B}|\ge 2$, e.g. for
$m=3$ and $B= \frac{\beta_0}2 + i -\frac{\sqrt 6}2$.

\begin{figure}[h]
\includegraphics[height=6cm]{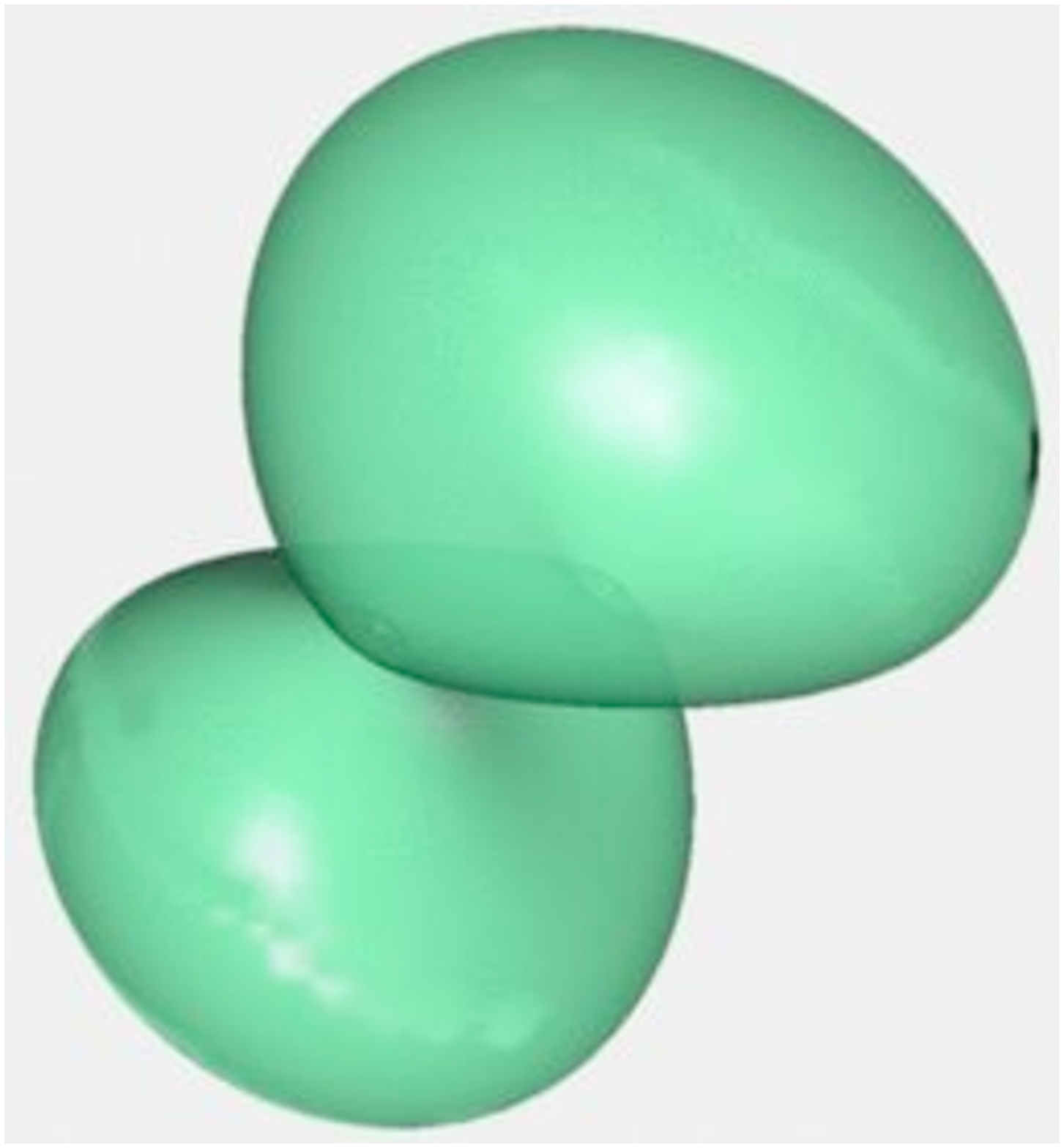}
\includegraphics[height=6cm]{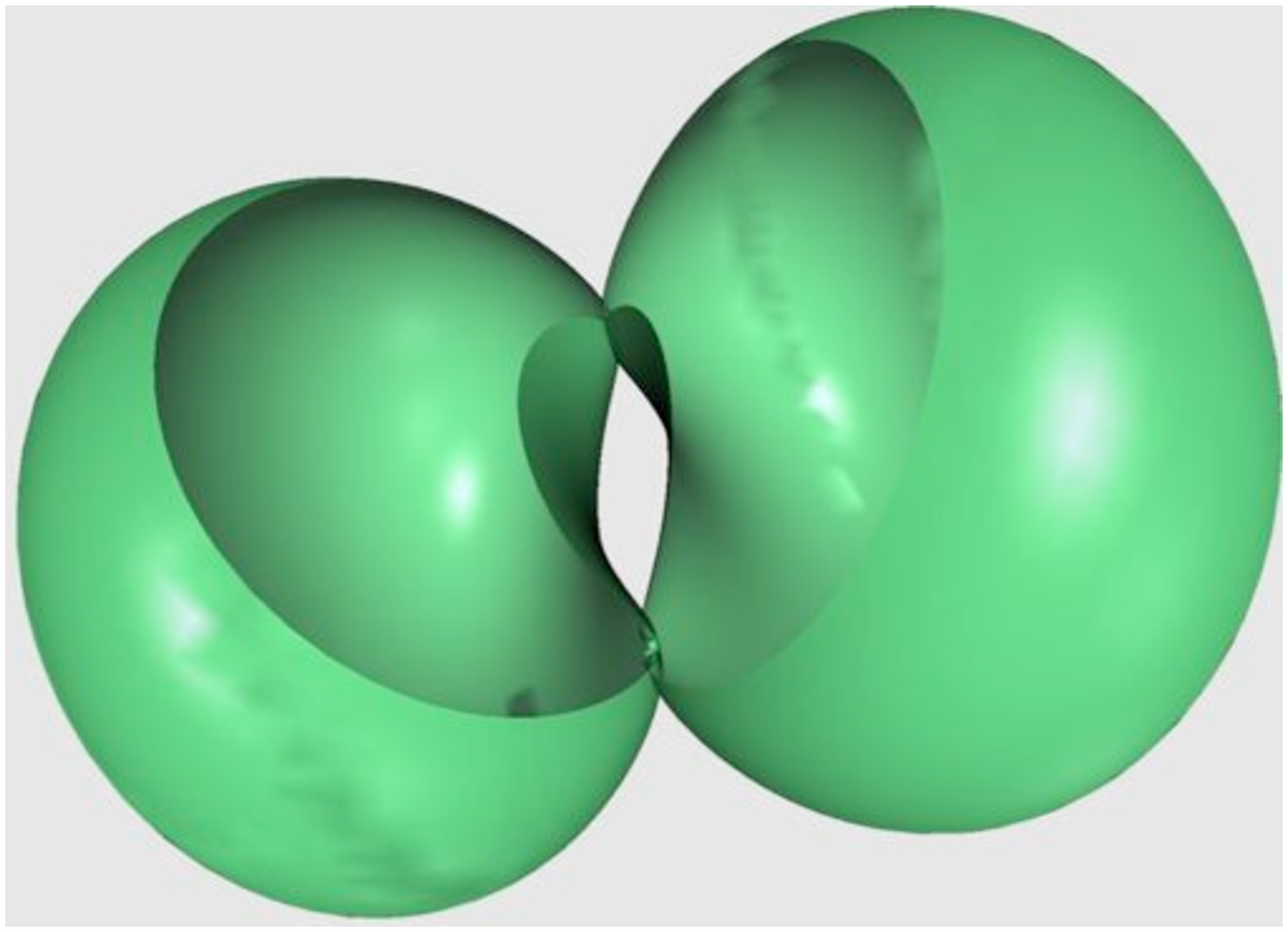}
\caption{Polychromatic Darboux transform of a homogeneous tori with
  $r_1=3$ and $r_2=1$, $B= \frac{\beta_0}2 + i -\frac{\sqrt 6}2$}
\end{figure}

\subsection{A Castro--Urbano torus}
Let us contemplate another example, which is part of a one-homogeneous
family found by Castro and Urbano \cite{castro_urbano} (see section
3.3 of \cite{pascal_frederic} for details). Take
$\Gamma=\mathrm{Span}\{1,i\}$ and $\beta_0=3+i$. Then
\[
\Gamma^*_{\beta_0} = \left\{ \frac{ 1\pm 3i }{2}, \frac{ -1\pm 3i }{2},\pm \frac{ 3 +i }{2} \right\}
\]
Using Theorem \ref{thm:hslLagrangianangle} we obtain a Hamiltonian
stationary torus by taking holomorphic sections without multiplier,
e.g., 
\[
f = \alpha_{\frac{-1+3i}2} + \alpha_{\frac{3+i}2}  = \frac 15e^{\pi j(3x-y)}\left((1 -4j + 3k)e^{2\pi i(3y-x)} + 
(1-3j +4k) e^{2\pi i(3x+y)}\right)
\]
is a 31--Castro--Urbano torus.  One easily verifies that $df= e^{\pi
  j(3x-y)} dz g$ with 
\[g= \frac{4\pi}5(3+i+3j-k)e^{\pi i (3y-x)} +
\frac{3\pi}5(3+i+j+3k)e^{\pi i(3x + y)}\,.
\]
\begin{figure}[h]
\includegraphics[height=6cm]{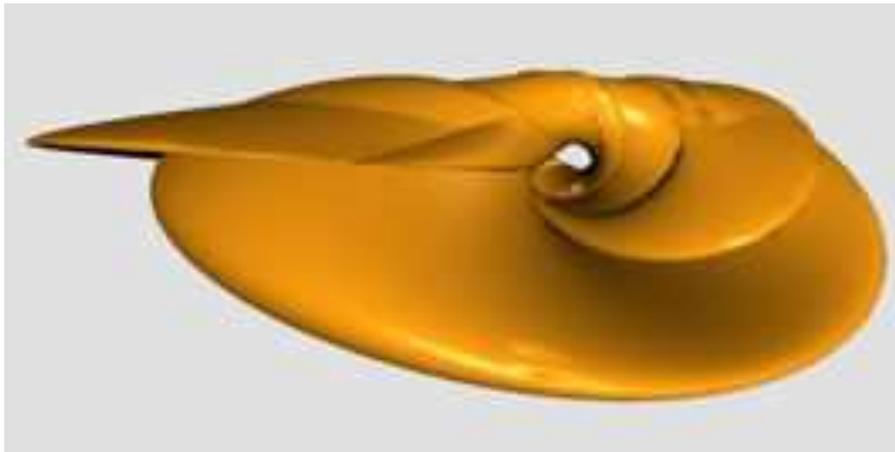}

\caption{Stereographic projection to $\R^3$ of a 31--Castro--Urbano torus 
}
\end{figure}
By Theorem \ref{cor:all monochromatic DT} all monochromatic Darboux
transforms of $f$ must be contained in the 31--Castro--Urbano family,
however in this case the Darboux transforms are not congruent to $f$.

\begin{figure}[h]
\includegraphics[height=6cm]{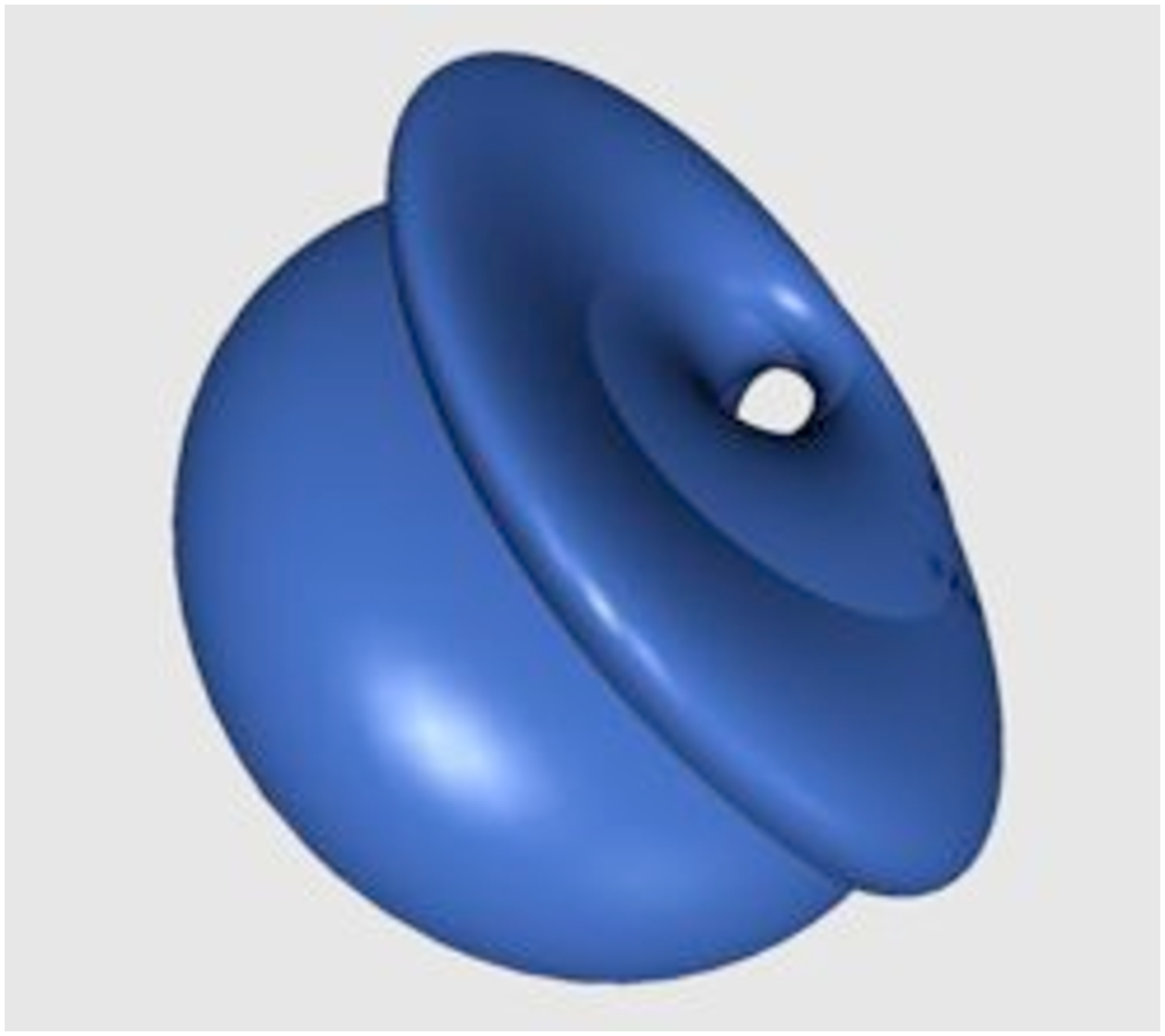}
\includegraphics[height=6cm]{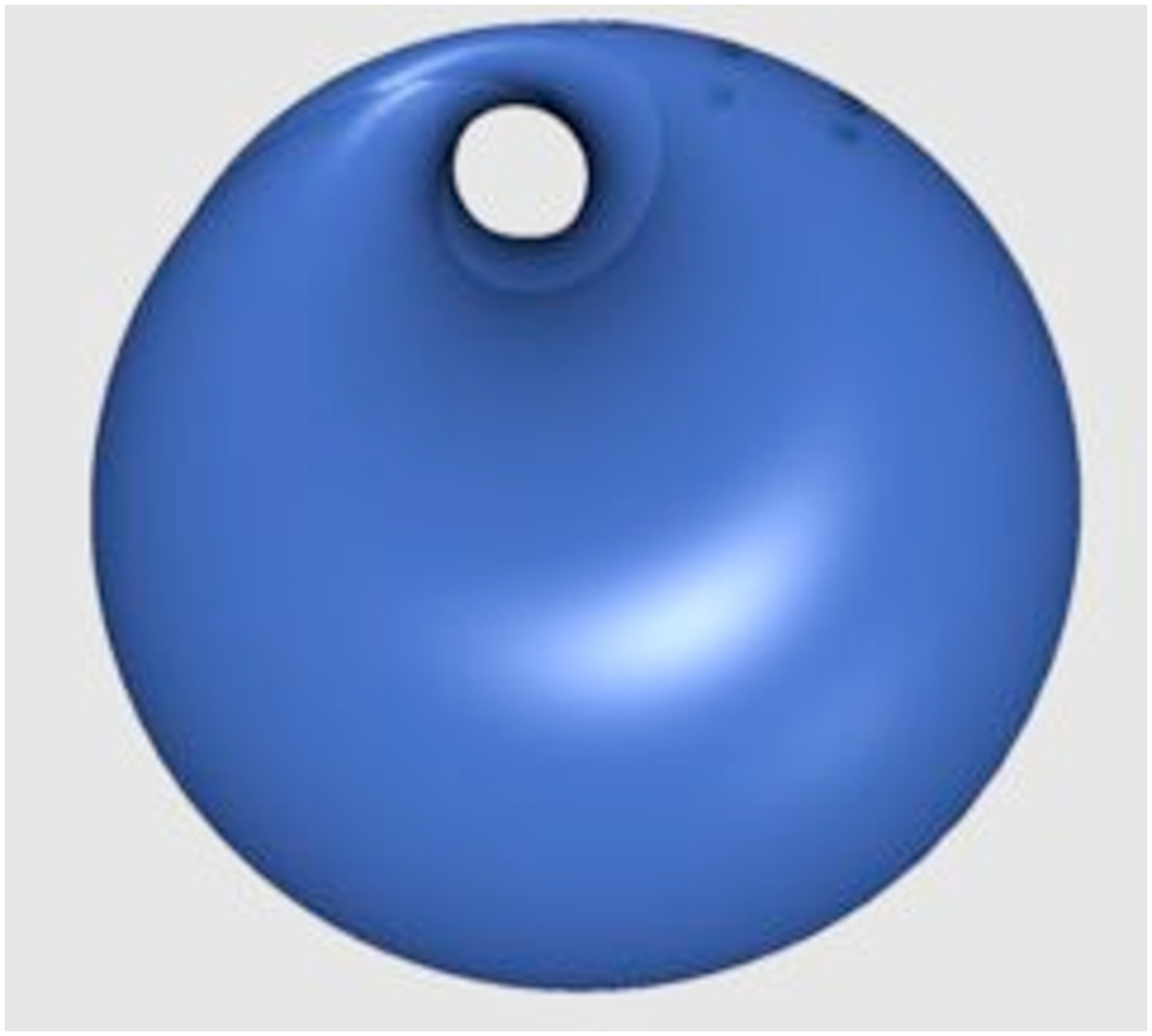}
\caption{Two monochromatic Darboux transforms of a 31--Castro--Urbano torus}
\end{figure} 

Again, there exists polychromatic Darboux transforms which are not
Lagrangian surfaces in $\Cj^2$.
\begin{figure}[h]
\includegraphics[height=6cm]{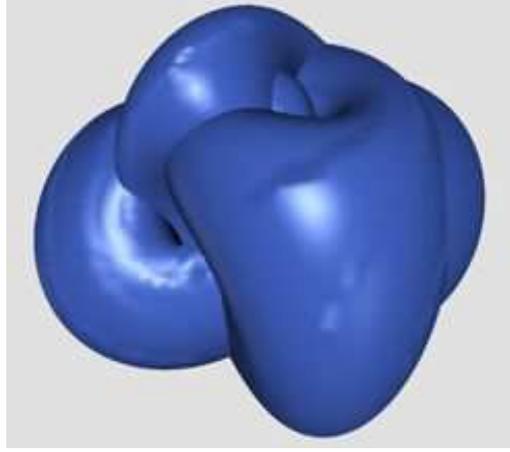}
\caption{Polychromatic Darboux transform of a 31--Castro--Urbano torus}
\end{figure}
Figure 9 shows a polychromatic Darboux transform of $f$ obtained when
choosing $B=\frac{\beta_0}2 + i -\frac{\sqrt{6}}2$ and
$\delta_1=\frac{\beta_0}2,
\delta_2=\frac{\beta_0}2+2i\in\Gamma^*_{0,B}$.

%%% Local Variables: 
%%% mode: latex
%%% TeX-master: "doc"
%%% End: 

%% file: mudarboux.tex
\section{$\mu$--Darboux transforms}
\label{sec:mu Darboux}

The Darboux transformation is defined for all conformal (branched)
immersions $f: M \to S^4$ of a Riemann surface $M$ into the
4-sphere. Whereas our previous discussion restricted to the case when
$M=T^2$ is a 2--torus we will now turn to the general Darboux
transformation of a Hamiltonian stationary immersion $f: M\to \R^4$ of
a Riemann surface $M$ into 4--space. In particular, we investigate the
link between Darboux transforms and the associated family of complex
flat connections given by the harmonic left normal $N =e^{j\beta}i$ of
$f$.  Given a Lagrangian immersion $f: M \to\R^4$ with Lagrangian
angle $\beta: M \to\R$ and $df = \ejbh dz g$, the left normal of $f$
is given by $N=e^{j\beta}i$. Since
\begin{equation}\label{eq:dnprime}
(dN)' = \frac{1}{2}(dN - N*dN) =
\frac j2 \ejbh(id\beta-*d\beta) \emjbh
\end{equation}
and
\[
 (dN)'' =  \frac{1}{2}(dN + N*dN) 
= \frac j2 \ejbh(id\beta+*d\beta) \emjbh
\]
we see that the real function $\beta$ is harmonic if and only if
$d(dN)'=0$, or, equivalently, $d(dN)''=0$, that is if $N$ is harmonic.
The mean curvature vector $\Hh$ of the immersion $f$ is given
\cite[Section 7.2]{coimbra} by 
\begin{equation}
\label{eq:H}
\Hh = N\bar H = \bar H R
\end{equation}
where $ (dN)' = -df H$. 
But (\ref{eq:dnprime}) simplifies with $\beta=2\pi \langle \beta_0, \cdot \rangle$ to 
$
(dN)' =  - \pi \ejbh dz \bar\beta_0 \ejbh k
$
and we see that $ H = \pi g\invers \bar\beta_0 \ejbh k\,.$

As customary \cite{uhlenbeck}, \cite{hitchin}, we can introduce a
spectral parameter $\mu\in\C_*$ so that the harmonicity of $N$ is
expressed in terms of the flatness of the family of complex
connections $d^\mu$. In our formulation, we consider the trivial
$\H$--bundle as a complex $\C^2$-bundle $(M \times \H, I) =
M\times\C^2$ where $I$ is the complex structure on $M\times \H$ given
by right multiplication with $i$. Then the $\C_*$--family of flat
complex connections is given \cite{Klassiker} by
\begin{equation}
\label{eq:flat connections}
d^\mu  = d + \frac 12 df H(N(a-1) + b)
\end{equation}
on $M\times\H$ where $\mu\in \C_*\subset \Span\{1, I\}$ and $a
=\frac{\mu+\mu\invers}{2}$, $b=\frac{\mu\invers-\mu }{2}I$.  Note that
by definition $a^2 + b^2=1$. Formally, the family of flat connections
is the same as in Moriya's paper \cite{katsuhiro}, however our
immersion is Hamiltonian in $\Cj^2$ rather than in $\C^2$ since we
choose a different complex structure on $\R^4$.

Furthermore, (\ref{eq:holomorphicity}) shows that every parallel
section $\alpha$ of $d^\mu$ is holomorphic since $d\alpha = -\frac 12
df H(N\alpha(a-1) + \alpha b)$.
\begin{definition}
  A branched conformal immersion $\hat f: M \to S^4$ which is given by
  the prolongation of a $d^\mu$--parallel section
  $\alpha\in\Gamma(\ttrivial)$ for some $\mu\in\C_*$ is called a
  \emph{$\mu$--Darboux transform} of a Hamiltonian stationary
  immersion $f: M \to\R^4$.
\end{definition}
Since a $\mu$--Darboux transform $\hat f= f+T$ is given by a
$d^\mu$--parallel section $\alpha$, Lemma \ref{lem:prolongation} shows
that $ T = \check TH\invers $ with
\begin{equation}
\label{eq: Tinvers}
\check T\invers =\frac{1}{2}( N\alpha(a-1)\alpha\invers + \alpha
b\alpha\invers)\,,
\end{equation}
We abbreviate $\check z = \alpha z\alpha\invers $ for
$z\in\C=\Span\{1,i\}$ so that $ T\invers = \frac{1}{2} H(N(\check a-1)
+ \check b) \,.  $ We now extend Theorem \ref{thm:Darboux is hsl} to
(local) $\mu$--Darboux transforms of a Hamiltonian stationary surface
$f: M \to\R^4$ where now $M$ may be an arbitrary Riemann surface.
\begin{theorem}
  Let $f:M \to\R^4$ be a Hamiltonian stationary immersion from a
  Riemann surface $M$ into the 4--space. Then the left normal of every
  (local) $\mu$--Darboux transform of $f$ is harmonic. In particular,
  every $\mu$--Darboux transform of $f$ is constrained Willmore.
\end{theorem}
 
\begin{proof} We essentially follow the proof in \cite{cmc_tori} that
  the Gauss map of a $\mu$--Darboux transform of a constant mean
  curvature surface is harmonic: From (\ref{eq:left normal darboux})
  we see that the left normal $\hat N$ of a Darboux transform $\hat f
  =f + T$ of $f$ is given by $\hat N = - T R T\invers$ where $R$ is
  the right normal of $f$.  Now (\ref{eq:H}) shows $R H = HN$, and we
  obtain $ \hat N = -\check T N \check T\invers\,, $ where $\check T$
  is given by a $d^\mu$--parallel section $\alpha$ and (\ref{eq:
    Tinvers}).  We put $\check \nu = \frac{1}{2}(N\alpha(a-1)+\alpha
  b)$ and compute its derivative, using $a^2 + b^2=1$, $dN = (dN)'' -
  dfH$ and $d\alpha = -df H \check \nu$, as
\[
d\check \nu = (dN)''\alpha\frac{a-1}2\,.
\]
Therefore, we obtain for $\check T\invers = \check \nu \alpha\invers$
the Riccati type equation
\begin{equation}
\label{eq:ricc}
d \check T\invers =
(dN)''\frac{\check a-1}2 + \check T\invers df H 
\check T\invers\,,
\end{equation}
and the derivative of $\hat N = - \check T N \check T\invers$ can  
be computed as
\[
d\hat N = (\hat N - N) df H \check T\invers + \check T(dN)''\left( N \frac{\check a-1}2 -\frac{\check a-1}2 \hat N \right) - \check T dN \check T\invers\,.
\]
The $(1,0)$--part of $d\hat N$ with respect to the complex structure
$\hat N$ is
\begin{equation}
\label{eq:dhatnprime}
(d\hat N)' = -\check T(dN)''\left(\frac{\check a-1}2 \hat N  + \frac{\check b}2\right) = -(dT + dfH)\,,
\end{equation}
where we used that $(\check a-1) \hat N + \check b = (1-\check a)
\check T$ by the definition (\ref{eq: Tinvers}) of $\check T$, and the
Riccati type equation (\ref{eq:ricc}).  In particular, since $N$ is
harmonic and $dfH=-(dN)'$ we see from \eqref{eq:dhatnprime} that
$(d\hat N)'$ is closed. In other words, the left normal $\hat N$ of a
$\mu$--Darboux transform of $f$ is harmonic.

We conclude the proof by showing that every conformal immersion
$f^\sharp: M \to \R^4$ with harmonic left normal $N^\sharp$ is
constrained Willmore.  From \cite{constrained} and the explicit
formulae for the mean curvature sphere congruence of an immersion in
\cite[Prop. 15]{coimbra} we see that $f^\sharp$ is constrained
Willmore if and only if there exists $\eta\in\Omega^1(\H)$ with
\begin{equation}
\label{eq:constr}
*\eta = -  R^\sharp\eta = \eta  N^\sharp\quad \text{ and } \quad d\left(\eta + d H^\sharp + R^\sharp*d H^\sharp + \frac 12 H^\sharp(N^\sharp dN^\sharp -*dN^\sharp)\right)=0\,,
\end{equation}
where $N^\sharp $ and $R^\sharp$ are the left and right normal of
$f^\sharp$, and $H^\sharp$ is given by $(dN^\sharp)'=\frac
12(dN^\sharp - N^\sharp *dN^\sharp)=-df^\sharp H^\sharp$. Since
$N^\sharp$ is harmonic we have $-df^\sharp \wedge dH^\sharp = d(d
f^\sharp H^\sharp)=0$ which implies that $*dH^\sharp = -R^\sharp
dH^\sharp$.  Therefore, $\eta\assign -\frac 12H^\sharp(N^\sharp
dN^\sharp -*dN^\sharp)$ is a 1--form satisfying (\ref{eq:constr}), and
every branched conformal immersion with harmonic left normal is
constrained Willmore. In particular, every $\mu$--Darboux transform of
a Hamiltonian Lagrangian immersion $f$ is constrained Willmore.
\end{proof}

In the case when $M=T^2$ is a 2--torus, the closed $\mu$--Darboux
transforms of a Hamiltonian stationary torus are exactly the
monochromatic Darboux transforms: To obtain parallel sections of
$d^\mu$ we use again the gauge $ \alpha= \ejbh \tilde \alpha$, and perform
a similar computation as in \cite{katsuhiro}.  From $d^\mu\alpha=0$ we
see $d\alpha = (dN)'\check T\alpha$ and (\ref{eq:dnprime}) then shows
\[
d\tilde\alpha = - \frac j4\big( 2d\beta\tilde \alpha 
+ (d\beta + i*d\beta)\tilde \alpha(a-1) + (*d\beta - i d\beta) \tilde \alpha b\big)\,.
\]
Decomposing $\tilde\alpha= \tilde\alpha_0+j\tilde\alpha_1$ we get the
 differential equation
\[
d\begin{pmatrix}\tilde \alpha_0\\
\tilde \alpha_1
\end{pmatrix}
= \frac 14 \begin{pmatrix} 0 & d\beta(\mu+1) - *d\beta i(\mu-1)\\
- \left(d\beta (\mu\invers + 1) + *d\beta i(\mu\invers-1)\right)&0
\end{pmatrix}\begin{pmatrix}\tilde \alpha_0\\
\tilde \alpha_1
\end{pmatrix}\,,
\]

where we used $\mu = a + ib$ and $\mu\invers = a- ib$.  If $M=T^2$ is
a 2--torus, the solution space of the differential equation is for
each $\mu\in\C_*$ 2--dimensional: writing $\beta(z) =
2\pi\langle\beta_0,z\rangle$ the general solution is
\begin{equation*}
\tilde \alpha =\begin{pmatrix}\tilde \alpha_0\\
\tilde \alpha_1
\end{pmatrix} =
\begin{pmatrix} 1\\ i \sqrt{\mu}\invers
\end{pmatrix}c_+e^{2\pi(x u + y v) }
+\begin{pmatrix} 1\\ -i \sqrt{\mu}\invers
\end{pmatrix} c_-e^{-2\pi(x u+ y v) } 
\end{equation*}
with $c_\pm\in\C$,
$
u =  \frac {i}{4} (\bar\beta_0 \sqrt\mu + \beta_0\sqrt{\mu}\invers)
$
and
$
v= \frac 1 4(\beta_0 \sqrt{\mu}\invers - \bar \beta_0 \sqrt{\mu})\,.
$
Putting 
\begin{equation}
\label{eq: A_mu}
 A^\mu=
\Re u + i \Re v = \frac{i\beta_0}4(\sqrt{\mu}\invers - \overline{\sqrt{\mu}})
\end{equation}
and
 \begin{equation}
\label{eq: C_mu}
 C^\mu = \Im u+ i
\Im v= \frac{\beta_0}4(\sqrt{\mu}\invers + \overline{\sqrt{\mu}})\,,
\end{equation}
we see that
\begin{equation}
\label{eq: eigenlines}
\alpha^\mu_\pm(z) = e^{\frac{j\beta(z)}{2}}(1\mp k \sqrt{\mu}\invers)
e^{\pm2\pi(\langle A^\mu,    z\rangle + i \langle C^\mu,z\rangle)}
\end{equation}
are holomorphic sections with multiplier 
\begin{equation}
\label{eq:multiplier}
h^\mu_\pm(\gamma) = e^{\pm 2\pi(\langle A^\mu,\gamma\rangle + i\langle\mp \frac{\beta_0}2 + C^\mu,\gamma\rangle)}\,.
\end{equation}
In other words, every $\mu\in\C_*$ gives two monochromatic, $\C$--independent
holomorphic sections $\alpha^\mu_\pm$ with multipliers $h^\mu_\pm$ satisfying
(\ref{eq:monochromatic holomorphic section}) with
$\delta=\frac{\beta_0}2$, $A=\pm A^\mu$, $\delta- B= \pm
C^\mu$, and
\begin{equation}
\label{eq:spec_parameter}
\lambda_\delta = \frac{2}{\beta_0}(\delta-iA - B) = \pm \frac 1{\sqrt{\mu}}\,.
\end{equation}
Conversely, we see that every monochromatic holomorphic
section is parallel for some $d^\mu$:

\begin{lemma}
\label{lem: spectral = mu}
For every monochromatic holomorphic section $\alpha$ there is a unique
$\mu\in\C_*$ such that $\alpha$ is $d^\mu$--parallel.
\end{lemma}
\begin{proof}
  Every monochromatic holomorphic section $\alpha\in H^0_h(\ttrivial)$
  with multiplier $h$ is given (\ref{eq:monochromatic holomorphic
    section}) by a complex scale of
\[
 \alpha_\delta = \ejbh(1 - k \lambda_\delta)e_{\delta-B} e^{2\pi \langle A, \cdot \rangle}
\]
where $\delta \in\Gamma^*_{A,B}$, $\lambda_\delta=
\frac{2}{\beta_0}(\delta -i A -B)$ and $ h=h^{A,B}$.  Since $\delta-B$
is independent of the choice of pair $(A,B)$ with $h=h^{A,B}$ both
$\lambda_\delta=\frac{2}{\beta_0}(\delta-iA -B)$ and $\mu=
\frac{1}{\lambda_\delta^2}$ are uniquely defined by $\alpha$.  As in
the proof of Proposition \ref{prop:normalization} we obtain from
(\ref{eq:monochromatic holomorphic section}) and (\ref{eq:lambda
  delta}) that
 \[
 \alpha_\delta = \ejbh(1\mp k \sqrt{\mu}\invers)e_{\delta-B} e^{2\pi \langle A, \cdot \rangle}
\]
with $A = \pm \frac{i\beta_0}4(\sqrt{\mu}\invers -
\overline{\sqrt{\mu}})$ and $\delta-B= \pm
\frac{\beta_0}4(\sqrt{\mu}\invers +
\overline{\sqrt{\mu}})$. Comparing with (\ref{eq: eigenlines}) we see
that $\alpha_\delta=\alpha^\mu_+$ or $\alpha_\delta =
\alpha^\mu_-$. But this proves that every monochromatic holomorphic
section $\alpha=\alpha_\delta c$, $c\in\C$, is parallel with respect
to $d^\mu$.
\end{proof}

\begin{cor}
\label{cor:2dimensionalspace}
The monodromy $H^\mu$ of the flat connection $d^\mu$ is a complex
multiple of the identity if and only if $\mu\in S^1$ with
$\beta_0\overline{\sqrt{\mu}}\in\Gamma^*$.
\end{cor}
\begin{proof}
  For a given $\mu$ we have a 2--dimensional space of parallel
  sections with the same multiplier if and only if $h^\mu_+ = h^\mu_-$
  which is equivalent with (\ref{eq:multiplier}), (\ref{eq: A_mu}) and
  (\ref{eq: C_mu}) to $A^\mu=0$ and $C^\mu =
  \frac{\beta_0}2\overline{\sqrt{\mu}}\in \frac 12 \Gamma^*$. Finally,
  $A^\mu=0$ is equivalent to $\mu\in S^1$.
\end{proof}
From (\ref{eq:flat connections}) we see that for $\mu\in S^1$ the
connection $d^\mu$ is in fact quaternionic. In particular, parallel
sections are quaternionically dependent. In other words, for all
$\mu\in S^1$ every $d^\mu$--parallel section gives the same
$\mu$--Darboux transform $\hat f$. In particular, even if $H^\mu$ is a
multiple of the identity, we can choose either one of $\alpha^\mu_\pm$
to get the $\mu$--Darboux transform $\hat f$ by prolongation.  We
summarize the previous discussions:
\begin{theorem}
\label{thm: spectral = mu}
Every closed $\mu$--Darboux transform $\hat f: T^2\to\R^4$ of a
Hamiltonian stationary torus $f: T^2\to\R^4$ is a monochromatic
Darboux transform of $f$, and vice versa.
\end{theorem}
\begin{rem} Again, we should point out that we only consider Darboux transforms closing on the original lattice.
\end{rem}

%%% Local Variables: 
%%% mode: latex
%%% TeX-master: "doc"
%%% End: 

%% file: spectralcurve.tex
\section{The spectral curve}

In this section we show that the multiplier spectral curve $\Sigma$ of
a Hamiltonian stationary torus $f: T^2\to \R^4$ and the spectral curve
$\Sigma_e$ of the harmonic complex structure given by $f$ coincide:
Given the $\C_*$--family of complex flat connections $d^\mu$ of a
harmonic complex structure $J$, the spectral curve of $J$ is defined
to be the normalization $\Sigma_e$ of the set of eigenvalues
\[
\Eig=\{(\mu, h^\mu) \mid \text{ there exists } \alpha\in\Gamma(\ttrivial) 
\text{ with } d^\mu\alpha =0, \gamma^*\alpha=\alpha h^\mu_\gamma\}\,.
\]
of the monodromy $H^\mu$ of the flat connections $d^\mu$. The
eigenlines of $H^\mu$ over points where $H^\mu$ is diagonalizable with
two different eigenvalues $h^\mu_\pm$ extend \cite{hitchin} to a line
bundle $\E$ over $\Sigma_e$.  In our situation, when $J$ is the
complex structure induced by the harmonic left normal of a Hamiltonian
stationary torus, the spectral curve of $J$ has spectral genus zero
\cite{katsuhiro}:

\begin{theorem}
  Let $f: T^2\to\R^4$ be a Hamiltonian stationary torus and $J$ the
  harmonic complex structure on $\trivial{}$ given by left
  multiplication by the left normal $N=e^{j\beta}i$ of $f$. Then the
  spectral curve $\Sigma_e$ of $J$ compactifies to
  $\bar\Sigma_e=\CP^1$ by adding points $x_0$ and $x_\infty$ over
  $\mu=0$ and $\mu=\infty$ respectively. Moreover, the map $\mu:
  \bar\Sigma_e\to \CP^1, (\mu,h^\mu)\mapsto \mu$ is a 2--fold covering
  over $\CP^1$ branched over $0$ and $\infty$, and the eigenline
  bundle $\E$ extends holomorphically to $\bar \Sigma_e$.
\end{theorem}
\begin{proof} 
We denote by
\[
\Eig_0=\{ (\mu, h^\mu) \mid H^\mu \not= h^\mu\Id\}
\]
the subset of $\Eig$ such that every $x=(\mu, h^\mu)\in\Eig_0$ has a
unique, up to complex scale, $d^\mu$--parallel section with multiplier
$h^\mu$.  We already have seen in Corollary
\ref{cor:2dimensionalspace} that $\Eig_0 =\Eig\setminus \Eig^0$,
where
\[
\Eig^0=\{(\mu,h^\mu) \mid \mu\in S^1 \text{ with }
\beta_0\overline{\sqrt{\mu}}\in\Gamma^*\}
\]
is finite.  By (\ref{eq: eigenlines}) the eigenlines $\E_{(\mu,
  h^\mu_\pm)}= \alpha^\mu_\pm\C$, $(\mu, h^\mu)\in\Eig_0$, extend
holomorphically to $\Eig^0$, and thus define a line bundle $\E$ over
$\Eig$.  Furthermore, for $\mu\in S^1$ with
$\beta_0\sqrt{\bar\mu}\in\Gamma^*$ we have (\ref{eq: eigenlines}) two
distinct limiting lines which shows that the normalization of $\Eig$
are two distinct copies of $\C_*$.  Since
\[
\E_{(\mu, h^\mu_\pm)} = \alpha_\pm \C= \ejbh (1\mp
k\sqrt{\mu}\invers)\C= \ejbh (\mp i\sqrt{\mu} + j)\C\,,
\]
and $\mu(x)\to 0$ for $x\to x_0$ and $\mu(x) \to \infty$ for $x\to
x_\infty$ respectively, we see that
\[
\E_{(\mu, h^\mu_\pm)} \to \E_{x_0}\assign\ejbh j\C, \quad (\mu, h^\mu_\pm)\to x_0\,,
\]
and
\begin{equation}
\label{eq:lineoverinfinty}
\E_{(\mu, h^\mu_\pm)} \to \E_{x_\infty} \assign \ejbh \C, \quad (\mu, h^\mu_\pm) \to x_\infty\,.
\end{equation}
\end{proof}

\begin{rem}
  There is a fixed point free real structure $\rho(\mu,h^\mu) \assign
  (\frac{1}{\bar\mu}, \bar h^\mu)$ on $\bar\Sigma_e$ which is
  compatible with the eigenline bundle $\E$ since
  $d^{\frac{1}{\bar\mu}}(\alpha j) = (d^\mu\alpha) j$.
\end{rem}

\begin{lemma}
  Let $\Sigma_e$ be the spectral curve of the harmonic complex
  structure $J$ of a Hamiltonian stationary torus $f$, and $\E$ the
  eigenline bundle over $\Sigma_e$.  Then the complex structure $J$
  can be recovered from $\bar\Sigma_e$ and $\E$ by quaternionically
  extending $J \varphi = \varphi i$,
  $\varphi\in\E_{x_\infty}$. Moreover, the Hamiltonian stationary
  torus $f$ is obtained as the limit of the Darboux transforms $f^x$
  when  $x\in\Sigma_e$ goes to $ x_\infty=\mu\invers(\infty)\in\bar\Sigma_e$.
\end{lemma}
\begin{proof}The first statement follows from $ J\ejbh = N \ejbh =
  \ejbh i $ and (\ref{eq:lineoverinfinty}). For the second, we observe
  that (\ref{eq: A_mu}) shows that $|A^\mu| =
  \frac{|\beta_0|}2|\sqrt{\mu}\invers-\overline{\sqrt{\mu}}|\to
  \infty$ as $\mu\to\infty$ and (\ref{eq:
    monochromatic tau, else}) therefore gives $f^x = f + T^x$ with
\[
|T^x| = \frac{|g|}{\pi\sqrt{4|A^\mu|^2+ \langle\beta_0,\frac{A^\mu}{|A^\mu|} \rangle^2}} \to 0, \quad  x\to x_\infty, \quad \mu=\mu(x)\,.
\]
\end{proof}

Since parallel sections are
  holomorphic we have a natural holomorphic map
\[
h\colon  \Eig \to \Spec, (\mu, h^\mu) \mapsto h^\mu\,.
\]
Lemma \ref{lem: spectral = mu} and Theorem \ref{thm:parametrized
  multipliers} show that this map is surjective.  By definition, every
point $x\in \Eig_0$ gives, up to complex scale, a unique monochromatic
holomorphic section with multiplier $h=h(x)$. If $x, \tilde
x\in\Eig_0$ with $h(x) = h(\tilde x)$ then $\dim_\C
H^0_h(\ttrivial)>1$. On the other hand, $H^0_h(\ttrivial)$ is spanned
by monochromatic holomorphic sections, and each monochromatic
holomorphic section $\alpha$ defines by Lemma \ref{lem: spectral =
  mu} a unique $\mu$ with $d^\mu\alpha=0$.  In particular, for
$h\not\in\Spec_0$, that is $\dim_\C H^0_h(\ttrivial)\not=1$, the
preimage of $h$ in $\Eig_0$ is finite.  Since $\Spec\setminus\Spec_0$
is a discrete subset in $\Spec$ this shows that $h|_{\Eig_0}$ is
injective away from a discrete set, and thus $h: \Eig\to\Spec$ is also
injective away from a discrete set since $\Eig^0=\Eig\setminus\Eig_0$
is finite.  Therefore, the induced map $h: \Sigma_e \to \Sigma$ on the
normalizations of $\Eig$ and $\Spec$ respectively is biholomorphic.
 
\begin{theorem}
  The multiplier spectral curve $\Sigma$ of a Hamiltonian stationary
  torus $f:T^2\to\R^4$ and the eigenline spectral curve $\Sigma_e$ of
  the harmonic left normal $N$ of $f$ are biholomorphic. In
  particular, the multiplier spectral curve compactifies to
  $\bar\Sigma=\CP^1$.

  Moreover, the kernel bundle $\L$ of $f$ and the eigenline bundle
  $\E$ of $N$ coincide, and the spectral curve $\Sigma$ parametrizes a
  $\CP^1$--family of Hamiltonian stationary tori isospectral to $f$.
\end{theorem}

%%% Local Variables: 
%%% mode: latex
%%% TeX-master: "doc"
%%% End: 

%% file: HSL_tori_final.bbl
\providecommand{\bysame}{\leavevmode\hbox to3em{\hrulefill}\thinspace}
\providecommand{\MR}{\relax\ifhmode\unskip\space\fi MR }
\providecommand{\MRhref}[2]{%
  \href{http://www.ams.org/mathscinet-getitem?mr=#1}{#2}
}
\providecommand{\href}[2]{#2}
\begin{thebibliography}{10}

\bibitem{anciaux}
H.~Anciaux, \emph{Construction of many {H}amiltonian stationary {L}agrangian
  surfaces in {E}uclidean four-space}, Calc. Var. Partial Diff. Equations
  \textbf{17} (2003), no.~2, 105--120.

\bibitem{holly}
H.~Bernstein, \emph{Non-special, non-canal isothermic tori with spherical lines
  of curvature}, Trans. Am. Math. Soc \textbf{353} (2001), 2245--2274.

\bibitem{conformal_tori}
C.~Bohle, K.~Leschke, F.~Pedit, and U.~Pinkall, \emph{Conformal maps from a
  2--torus to the 4--sphere}, arXiv:0712.2311.

\bibitem{constrained}
C.~Bohle, G.~Peters, and U.~Pinkall, \emph{Constrained {W}illmore surfaces},
  Calculus of Variations and Partial Differential Equations \textbf{32} (2008),
  no.~2, 263--277.

\bibitem{coimbra}
F.~Burstall, D.~Ferus, K.~Leschke, F.~Pedit, and U.~Pinkall, \emph{Conformal
  geometry of surfaces in ${S}^4$ and quaternions}, Lecture Notes in
  Mathematics, Springer, Berlin, Heidelberg, 2002.

\bibitem{cmc_tori}
E.~Carberry, K.~Leschke, and F.~Pedit, \emph{Darboux transforms and spectral
  curves of constant mean curvature surfaces revisited}, in preparation.

\bibitem{castro_chen}
I.~Castro and B.-Y. Chen, \emph{Lagrangian surfaces in complex {E}uclidean
  plane via spherical and hyperbolic curves}, Tohoku Math. J, 2nd ser
  \textbf{58} (2006), no.~4, 565--579.

\bibitem{castro_urbano}
I.~Castro and F.~Urbano, \emph{Examples of unstable {H}amiltonian-minimal
  {L}agrangian tori in $\mathbb{C}^2$}, Compositio Mathematica \textbf{111}
  (1998), 1--14.

\bibitem{chen}
B.-Y. Chen, \emph{Construction of {L}agrangian surfaces in complex {E}uclidean
  plane with {L}egendre curves}, Kodai Math. J. \textbf{29} (2006), no.~1,
  84--112.

\bibitem{chen_morvan}
B.-Y. Chen and J.-M. Morvan, \emph{Deformations of isotropic submanifolds in
  {K}\"ahler manifolds}, J. of Geometry and Physics \textbf{13} (1994),
  79--104.

\bibitem{darboux}
G.~Darboux, \emph{Sur les surfaces isothermiques}, C. R. Acad. Sci. Paris
  \textbf{128} (1899), 1299--1305.

\bibitem{Klassiker}
D.~Ferus, K.~Leschke, F.~Pedit, and U.~Pinkall, \emph{Quaternionic holomorphic
  geometry: {P}l\"ucker formula, {D}irac eigenvalue estimates and energy
  estimates of harmonic 2-tori}, Invent. math. \textbf{146} (2001), 507--593.

\bibitem{schmidt_grinevich}
P.~Grinevich and M.~Schmidt, \emph{Conformal invariant functionals of
  immersions of tori in $\mathbb{R}^3$}, J. Geom. Phys. \textbf{26} (1998),
  51--78.

\bibitem{pascal_quaternionic}
F.~H\'{e}lein and P.~Romon, \emph{Weierstrass representation of {L}agrangian
  surfaces in four dimensional space using spinors and quaternions},
  Commentarii Mathematici Helvetici \textbf{75} (2000), 668--680.

\bibitem{pascal_frederic}
\bysame, \emph{Hamiltonian stationary {L}agrangian surfaces in $\mathbb{C}^2$},
  Communications in Analysis and Geometry \textbf{10} (2002), no.~1, 79--126.

\bibitem{hitchin}
N.~Hitchin, \emph{Harmonic maps from a $2$-torus to the $3$-sphere}, J.
  Differential Geom. \textbf{31} (1990), no.~3, 627--710.

\bibitem{habil}
K.~Leschke, \emph{Transformation on {W}illmore surfaces}, Habilitationsschrift,
  Universit\"at Augsburg, 2006.

\bibitem{pascal_ian}
I.~McIntosh and P.~Romon, \emph{The spectral data for {H}amiltonian stationary
  {L}agrangian tori in $\mathbb{R}^4$}, arXiv:0707.1767.

\bibitem{katsuhiro}
K.~Moriya, \emph{Hamiltonian stationary {L}agrangian tori in the complex
  {E}uclidean plane with rational spectral curve}, arXiv:0710.4233.

\bibitem{oh}
Y.-G. Oh, \emph{Volume minimization of {L}agrangian submanifolds under
  {H}amiltonian deformations}, Math. Z. \textbf{212} (1993), 175--192.

\bibitem{taimanov_weierstrass}
I.~Taimanov, \emph{The {W}eierstrass representation of closed surfaces in
  $\mathbb{R}^3$}, Funct. Anal. Appl. \textbf{32} (1998), 49--62.

\bibitem{uhlenbeck}
K.~Uhlenbeck, \emph{{H}armonic maps into {L}ie groups (classical solutions of
  the chiral model)}, J. Diff.\ Geom. \textbf{30} (1989), 1--50.

\end{thebibliography}
